\documentclass{article}
\usepackage[utf8]{inputenc}
\usepackage{amsmath,amsthm, amssymb, amsfonts}
\usepackage{todonotes}
\usepackage{kantlipsum}
\usepackage{wrapfig}
\usepackage{paralist}
\usepackage{wasysym}
\usepackage{multicol}
\usepackage{float} 
\usepackage{verbatim}
\usepackage{fullpage}
\usepackage{mathtools}
\usepackage{diagbox}
\usepackage[pdfencoding=auto]{hyperref}
\usepackage{graphicx}
\usepackage{wrapfig}
\usepackage{mathrsfs}
\usepackage{enumerate}
\usepackage{dsfont}
\usepackage{units}

\usepackage[]{algorithm2e}
\usepackage{ltablex}
\usepackage{varwidth}

\usepackage{xcolor}

\allowdisplaybreaks[1]
\usepackage{sistyle} 
\newcommand\bSI[1]{{\small[\SI{}{#1}]}}

\makeatletter
\newlength\unitwdth
\newlength\numwdth
\newlength\tdima
\newcommand\SIdescr[2]{%
    \setlength\tdima{\linewidth}%
    \addtolength\tdima{\@totalleftmargin}%
    \addtolength\tdima{-\dimen\@curtab}%
    \addtolength\tdima{-\unitwdth}%
    \addtolength\tdima{-\numwdth}%
    \parbox[t]{\tdima}{%
        #1
        \leaders\hbox{$\m@th\mkern \@dotsep mu\hbox{\tiny.}\mkern \@dotsep mu$}%
        \hfill
        \ifhmode\strut\fi
        \makebox[0pt][l]{%
            \makebox[\unitwdth][l]{}%
            \makebox[\numwdth][r]{#2}}}}
\makeatother

\newcommand{\Ycal}{\mathcal{Y}}
\newcommand{\Hil}{\mathcal{H}}

\newcommand{\N}{\mathbb{N}}

\newcommand{\R}{\mathbb{R}}

\newcommand{\eps}{\varepsilon}
\newcommand{\Realization}{\mathrm{R}}

\newcommand{\Sol}{\mathcal{P}}

\newcommand\norm[1]{\left\lVert#1\right\rVert}
\newcommand\abs[1]{\left|#1\right|}

\newtheorem{theorem}{Theorem}[section]
\newtheorem*{theorem*}{Theorem}
\newtheorem{remark}[theorem]{Remark}
\newtheorem{Question}[theorem]{Question}
\newtheorem{definition}[theorem]{Definition}

\newtheorem*{remark*}{Remark}
\newtheorem*{proposition*}{Proposition}

\numberwithin{equation}{section}

\definecolor{darkcandyapplered}{rgb}{0.64, 0.0, 0.0}

\title{Numerical Solution of the Parametric Diffusion Equation \\ by Deep Neural Networks}

\author{Moritz Geist\footnote{These authors contributed equally.} \thanks{Institut für Mathematik, Technische Universit\"at Berlin, Straße des 17.~Juni 136, 10623 Berlin, Germany, e-mail: \texttt{$\{$geist, raslan,  schneidr, kutyniok$\}$@math.tu-berlin.de}} 
\qquad\quad  Philipp Petersen\footnotemark[1]  \thanks{University of Vienna, Faculty of Mathematics and Research Plattform Data Science @ Uni Vienna, Oskar Morgenstern Platz 1,
1090 Wien, e-mail: \texttt{philipp.petersen@univie.ac.at}} 
\qquad\quad Mones Raslan\footnotemark[1] \footnotemark[2] 
\\  Reinhold Schneider\footnotemark[2]  
\qquad\quad Gitta Kutyniok\footnotemark[2] \footnote{
Fakult\"at Elektrotechnik und Informatik, Technische Universit\"at Berlin}
\footnote{
Department of Physics and Technology, University of Troms\o 
}
}

\begin{document}

\maketitle
\begin{abstract}
    We perform a comprehensive numerical study of the effect of approximation-theoretical results for neural networks on practical learning problems in the context of numerical analysis. As the underlying model, we study the machine-learning-based solution of parametric partial differential equations. Here, approximation theory predicts that the performance of the model should depend only very mildly on the dimension of the parameter space and is determined by the intrinsic dimension of the solution manifold of the parametric partial differential equation. We use various methods to establish comparability between test-cases by minimizing the effect of the choice of test-cases on the optimization and sampling aspects of the learning problem. We find strong support for the hypothesis that approximation-theoretical effects heavily influence the practical behavior of learning problems in numerical analysis.
\end{abstract}



\section{Introduction}
This work studies the problem of numerically solving a specific parametric partial differential equation (PPDE) by training and applying neural networks (NNs). The central goal of the following exposition is to identify those key aspects of a parametric problem that render the problem harder or simpler to solve for methods based on NNs.

The underlying mathematical problem, the solution of PPDEs, is a standard problem in applied sciences and engineering. In this model, certain parts of a PDE such as the boundary conditions, the source terms, or the shape of the domain are controlled through a set of parameters, e.g., \cite{CertReduced, QuarteroniIntro}. In some applications where PDEs need to be evaluated very often or in real-time, individually solving the underlying PDEs for each choice of parameters becomes computationally infeasible. In this case, it is advisable to invoke methods that leverage on the joint structure of all the individual problems.
A typical approach is that of constructing a \emph{reduced basis} associated with the problem. With respect to this basis, the computational complexity of solving the PPDE is then significantly reduced, e.g., \cite{CertReduced, QuarteroniIntro, RozzaPateraRB, RedChal}. 

Recently, as an alternative or to augment the reduced basis method, approaches were introduced that attempt to learn the parameter-to-solution map through methods of machine learning. We will provide a comprehensive overview of related approaches in Section \ref{sec:RelWork}. One approach is to train a NN to fit the discretized parameter-to-solution map, i.e., a map taking a parameter to a finite-element discretization of the solution of the associated PDEs. This approach has already been analyzed theoretically in \cite{NNParametric} where it was shown from an approximation-theoretical point of view that the hardness of representing the parameter-to-solution map by NNs is determined by a highly problem-specific notion of complexity that depends (in some cases) only very mildly on the dimension of the parameter space. 

In this work, \emph{we study the problem of learning the discretized parameter-to-solution map in practice}. We hypothesize that the approximation-theoretical capacity of a NN architecture is one of the central factors in determining the difficulty level of the learning problem in practice.

The motivation for this analysis is two-fold: First, we regard this as a general analysis of the feasibility of approximation-theoretical arguments in the study of deep learning. Second, specifically for the problem of numerical solution of PPDEs, we consider it important to identify which characteristics of a parametric problem determine its practical hardness. This is especially relevant to identify in which areas the application of this model is appropriate. We outline these two points of motivation in Section \ref{sec:Motivation}. The design of the numerical experiment is presented in Section \ref{sec:ExperimentIntro} and we give a high-level report of our findings in Section \ref{sec:Findings}.

\subsection{Motivation}\label{sec:Motivation}
As already outlined before, we describe the motivation for this paper in the following two sections. 

\subsubsection{Understanding of Deep Learning in General}\label{sec:Mot1}

A typical learning problem consists of an unknown data model, a hypothesis class, and an optimization procedure to identify the best fit in the hypothesis class to the observed (sampled) data, e.g., \cite{Cucker02onthe, cucker_zhou_2007}.  In a deep learning problem, the hypothesis class is the set of NNs with a specific architecture.

The approximation-theoretical point of view analyzes the trade-off between the capacity of the hypothesis class and the complexity of the data model. In this sense, this point of view describes only one aspect of the learning problem.

In the framework of approximation theory, there are precise ways to assess the hardness of an underlying problem. Concretely, this is done by identifying the rate by which the misfit between the hypothesis class and the data model decreases for sequences of growing hypothesis classes. For example, one common theme in the literature is the observation that for certain function classes, NNs do not admit a curse of dimension, i.e., their approximation rates do not deteriorate exponentially fast with increasing input dimension, e.g., \cite{Barron1993, shaham2018provable, poggio2017and}. Another theme is that classes of smooth functions can be approximated more efficiently than classes of rougher functions, e.g., \cite{Mhaskar:1996:NNO:1362203.1362213, YAROTSKY2017103, PetV2018OptApproxReLU, FEMNNsPetersenSchwab}.

While these results offer some interpretation of why a certain problem should be harder or simpler, it is not clear how relevant these results are in practice. Indeed, there are at least three issues that call the approximation-theoretical explanation for a practical learning problem into question:

\begin{itemize}
    \item \textit{Tightness of the upper bounds:} Approximation-theoretical bounds usually describe worst-case error estimates for whole classes of functions. For individual functions or subsets of these function classes, there is no guarantee that one could not achieve a significantly better approximation rate. 
    
    \item \textit{Optimization and sampling prevent approximation theoretical effect from materializing:} As explained at the beginning of this section, the learning problem consists of multiple aspects, one of which is the ability of the hypothesis class to describe the data. 
    Two further aspects are how well the sampling of the data model describes the true model and how well the optimization procedure performs in finding the best fit to the sampled data. 
    Since the underlying optimization problem of deep learning is in general non-convex, it is conceivable that, while there theoretically exists a very good approximation of a function by a NN, finding it in practice is highly unlikely. Moreover, it is certainly possible that the sampling process does not contain sufficient information to guarantee that the optimization routine will identify the theoretically best approximation.
        
    \item \textit{Asymptotic estimates:} All approximation-theoretical results mentioned until here and almost all in the literature describe the capacity of NNs to represent functions approximately with accuracy $\eps$ for sufficiently large architectures only in a regime where $\eps$ tends to zero and the size of the architecture is sufficiently large. 
    The associated approximation rates may contain arbitrarily large implicit constants, and therefore it is entirely unclear if changes to the trade-off between the complexity of the data model and the size of the architecture have the theoretically predicted impact for moderately-sized practical learning problems.

\end{itemize}

We believe that, to understand the effect of approximation-theoretical capacities of NNs in practical learning scenarios, the learning problem associated with the parameter-to-solution map in a PPDE occupies a special role:
It is, in essence, a high-dimensional approximation problem of a function that has a very strong low-dimensional, but highly non-trivial structure. What is more is that one can, to a certain extent, control the complexity of the problem, as we have seen in \cite{NNParametric}.  
In this context, we can ask ourselves the following questions: Do we observe a curse of dimensionality in the practical solution of the problem? If not, how does the difficulty in practice scale with the parameter dimension? On which characteristics of the problem does the hardness of the practical solution thereof depend?

If we study these questions numerically, then the answers can be compared with the predictions from approximation-theoretical considerations. If the predictions coincide with the observed behavior and other causes, such as artefacts from the optimization and sampling procedure, can be ruled out, then we can view these experiments as a strong support for the practical relevance of approximation-theoretical arguments.

Because of this, we study the aforementioned questions in an extensive numerical experiment that will be described in Section \ref{sec:ExperimentIntro} below.

\subsubsection{Feasibility of the Machine-Learning-Based Solution of Parametric PDEs}\label{sec:Mot2}

The method (as described in \cite{NNParametric}) of learning the parameter-to-solution map has at least two major advantages over classical approaches to solve PPDEs: First of all, the setup is completely independent of the underlying PPDE. This versatility of NNs could be quite desirable in an environment where many substantially different PPDEs are treated. Second, because this approach is fully data-driven, we do not require any knowledge of the underlying PDE. Indeed, as long as sufficiently many data points are supplied, for example, from a physical experiment, the approach could be feasible under high uncertainty of the model. 

The main drawback of the method is the lack of theoretical understanding thereof. Moreover, for the theoretical results that do exist, we lack any evaluation of how pertinent the theoretical observations are for practical behavior. Most importantly, we do not have an a priori assessment for the practical feasibility of certain problems. 

In \cite{NNParametric}, we observed that the complexity of the solution manifold, i.e., the set of all solutions of the PDE, is a central quantity involved in upper bounding the hardness of approximating the parameter-to-solution map with a NN. In practice, it is unclear to what extent this notion is appropriate and if the complexity of the solution manifold influences the performance of the method at all.

In the numerical experiment described in the next chapter, we explore the performance of the learning approach for various test-cases with different intrinsic complexities and observe the sensitivity of the method to the different setups.

\subsection{The Experiment}
 \label{sec:ExperimentIntro}
To analyze the approximation-theoretical effect of the architecture on the overall performance of the learning problem in practice, we train a fully connected NN as considered in \cite{NNParametric} on a variety of datasets stemming from different parameter choices of the parametric diffusion equation. 
The design of the data sets is such that we vary the relationship between the capacity of the architecture and the complexity of the data and report the effect on the overall performance. 

In designing such an experiment, we face \emph{three fundamental challenges hindering the comparability between test-cases}:

\begin{itemize}
\item \textit{Effect of the optimization procedure:}
The effect of the architecture on the optimization procedure is not clear, and this interplay may be a much stronger factor in the performance of the method than the capacity of the architecture to fit the data model. 
Similarly, the effect of the complexity of the data model could affect the optimization procedure and influence the performance of the learning method stronger than any approximation-theoretical effect.

\item \textit{Effect of the sampling procedure:}
We train our network based on a finite number of samples of the true solution. The number and choice of samples could have a non-negligible effect on the overall performance and most importantly affect some test-cases more than others. 

\item \textit{Quantification of the intrinsic complexity:} 
While we have theoretically established that the complexity of the solution manifold is the main factor in upper-bounding the hardness of the problem in the approximation-theoretical framework, we cannot, in practice, quantify this complexity.
\end{itemize}

We address these issues in the following four ways:

\begin{itemize}
\item \textit{Keeping the architecture fixed:} 
An approximation-theoretical result on NNs is based on three ingredients. A function class $\mathcal{C}$, a worst-case accuracy $\eps>0$, and the size of the architecture. 

Whenever one of these hyper-parameters---the function class, the accuracy, or the architecture---is fixed, one can theoretically describe how changing a second parameter influences the last one. For example, for fixed $\mathcal{C}$, an approximation-theoretical statement yields an estimate of the necessary size of the architecture to achieve an accuracy of $\eps$.

Because of the potentially strong impact of the architecture on the optimization procedure, we expect that the most sensible point of view to test numerically is that where the architecture remains fixed while we vary the function class $\mathcal{C}$ and observe $\eps$. 
This way, we can guarantee that the influence of the architecture on the optimization procedure is the same between test-cases.

\item \textit{Analyzing the convergence behavior a posteriori:}
We are not aware of any method to guarantee a priori that the choice of the data model would not influence the convergence behavior. We do, however, analyze the convergence after the experiment to see if there are fundamental differences between our test-cases. This analysis reveals no significant differences between all the setups and therefore indicates that the effect of the data model on the optimization procedure is very similar between test-cases.

\item \textit{Establishing independence of sample generation:} We run the experiment multiple times for various numbers of training samples $N$ chosen in the same way---uniformly at random---in every test-case. Between the choices of $N$, we observe a linear dependence of the achieved accuracy on $N$. This indicates that the influence of the number of $N$ on the performance of the method is the same for all test-cases.

\item  \textit{Design of semi-ordered test-cases:} 
While we are not able to assess the intrinsic complexity exactly, it is straight-forward to construct series of test-cases with increasing complexity. In this sense, we can introduce a semi-ordering of test-cases according to their complexity and observe to what extent the performance of the method follows this ordering. 
\end{itemize}
We present the construction of the test-cases in Section \ref{sec:NumericalRes} and discuss the measures taken to remove effects caused by the optimization and sampling procedures in greater detail in Appendix \ref{sec:EliminationCauses}. All of our test-cases consider the following  parametric diffusion equation
\begin{align*}
	    -\nabla\cdot (a_y(\mathbf{x}) \cdot \nabla u_y(\mathbf{x})) = f(\mathbf{x}), \quad \text{ on } \Omega= (0,1)^2, \quad {u_y}|_{\partial  \Omega} = 0,
\end{align*}
where $f\in L^2(\Omega)$ and $a_y \in  L^{\infty}(\Omega)$, is a diffusion coefficient depending on a parameter $y \in \mathcal{Y}$.
In our test-cases below, we learn a discretization of the map $	\R^p \supset \mathcal{Y} \ni y \mapsto u_y$, 
where $p \in \N$, for various choices of parametrizations
\begin{align}\label{eq:ParamIntro}
	\R^p \supset \mathcal{Y} \ni y \mapsto a_y.
\end{align}
Concretely, we vary the following characteristics of the parametrizations and observe the effect on the overall performance of the learning problem:
\begin{itemize}
\item \textit{Type of parametrization:} We choose test-cases which differ with respect to the following characteristics: First, we study parametrizations \eqref{eq:ParamIntro} of various degrees of smoothness. Second, we study test-cases where the parametrization \eqref{eq:ParamIntro} is affine-linear and non-linear. Third, we consider cases, where $a_y = \sum_{i=1}^p \widetilde{a}_{y_i}$ for $\widetilde{a}_{y_i} \in L^{\infty}(\Omega)$ and where the supports of $(\widetilde{a}_{y_i})_{i=1}^p$ overlap or have various degrees of separation.
\item \textit{Dimension of parameter space:} 
The discretization of our solution space is done on the maximal computationally feasible grid (with respect to our workstation). We have chosen the dimensions $p$ of the parameter spaces in such a way that the resolutions of the parametrized solutions are still meaningful with respect to the underlying discretization. 

\item \textit{Complexity hyper-parameters:}
To generate comparable test-cases with increasing complexities, we include two types of hyper-parametes into the data-generation process. One that directly influences the ellipticity of the problem and another that introduces a weighting of the parameter values.
\end{itemize}

We expect that these tests yield answers to the following questions: How versatile is the approach? Does it perform well only for special types of parametrizations or is it generally applicable? Do we observe a curse of dimensionality and how much does the performance of the learning method depend on the dimension of the parameter space? How strongly does the performance of the learning method depend on the intrinsic complexity of the data?

\subsection{Our Findings}\label{sec:Findings}

In the numerical experiments, which we report in Section \ref{subsec:Results} and evaluate in Section \ref{subsec:Eval},  we find that the proposed method is very sensitive to the underlying type of test-case. Indeed, we observe qualitatively different scaling behaviors of the achieved error with the dimension $p$ of the parameter space  between different test-cases. Concretely, we observe the following asymptotic behavior of the errors in different test-cases: $\mathcal{O}(1), \mathcal{O}(\log(p))$ and $\mathcal{O}(p^k)$ for $p \to \infty$ and $k > 0$, where $k$ depends on one of the complexity hyper-parameters. Notably, we do not observe a scaling according to the curse of dimensionality, i.e., an error scaling exponentially with $p$, in any of the test-cases. 
We also observe that the achieved errors obey the semi-ordering of complexities of the test-cases. This shows that the method is very versatile and can be applied for various settings. Moreover, the complexity of the solution manifold appears to be a sensible predictor for the efficiency of the method.

In addition, we observe that the numerical results agree with the predictions that can be made via approximation-theoretical considerations. By design, we can exclude effects associated with the optimization and sampling procedures. This supports the practical relevance of approximation-theoretical results for this particular problem and for deep learning problems in general. 

\subsection{Related Works}\label{sec:RelWork}

The practical application of NNs in the context of PDEs dates back to the 1990s, \cite{PDE_NN_historic}. 
However, in recent years the topic again gained traction in the scientific community driven by the ever-increasing availability of computational power. 
Much of this research can be condensed into three main directions: 
Learning the solution of a single PDE, system identification, and goal-oriented approaches. 
The first of these directions uses NNs to directly model the solution of a (in some cases user-specified) single PDE, \cite{deepRitz,PhysicsInformed1, yang2018physics,DeepXDE, samaniego2019energy}, an SDE, \cite{BeckGrohsJentzen,weinan2017deep}, or even the joint solution for multiple boundary conditions, \cite{DGM}. 
These methods mostly rely on the differential operator of the PDE to evaluate the loss, but other approaches do exist, \cite{han2020derivativefree}. 
In system identification, one tries to discover an underlying physical law from data by reverse-engineering the PDE. 
This can be done by attempting to uncover a hidden parameter of a known equation, \cite{PhysicsInformed2}, or modeling physical relations, \cite{DDdiscovery, raissi2018deep}. 
Conversely, goal-oriented approaches, try to infer a quantity of interest stemming from the solution of an underlying PDE. 
For example, NNs can be used as a surrogate model to directly learn the quantity of interest and thereby circumvent the necessity of explicitly solving the equation, \cite{Khoo}. 
A practical example for this is given by the ground state energy of a molecule which is derived from the solution of the electronic Schr\"odinger equation. 
This task has been efficiently solved by graph NNs, \cite{ES1,ES2,ES3}. Furthermore, building a surrogate model can be especially useful in uncertainty quantification, \cite{DeepUQ}. 
NNs can also aid classical methods in solving goal-oriented tasks, \cite{brevis2020datadriven,MultiLevelML}. 
In addition to the aforementioned research directions, further work has been done on fusing NNs with classical numerical methods to assist, for example, in model-order reduction,  \cite{regazzoni2019, lee2018model}. 

Our work focuses on PPDEs and more specifically we are interested in learning the mapping from the parameter to the coefficients of the high-fidelity solution. Related but different approaches were analyzed in \cite{RBNonlinearProblems}, and \cite{DalSantoParametric, DeepUQ}, where the solution of the PPDE is learned in an already precomputed reduced basis or at point evaluations in fixed spatial coordinates.

On the theoretical side, the majority of works analyzing the power of NNs for the solution of (parametric) PDEs is concerned with an approximation-theoretical approach. 
Notable examples of such works include \cite{Curse, weinan2017deep, grohs2018proof, MachLearningSPDEJentzen, JentzenEHighSPDEs, SchwabOption, JentzenKolmogorov, BeckGrohsJentzen, berner2018analysis, JentzenHeat}, in which it is shown that NNs can overcome the curse of dimensionality in the approximative solution of some specific single PDE. 
In the same framework, it was shown in \cite{berner2018analysis} how estimates on the approximation error imply bounds on the generalization error.
Concerning the theoretical analysis of PPDEs, we mention \cite{SchwabZech,NNParametric,PetersenTransport,SchwabQoI}. 
We will describe the results of the first two works in more detail in Section \ref{sec:ApproxOfSolutionMapTheoretical}. The work \cite{SchwabQoI} is concerned with an efficient approximation of a map that takes a noisy solution of a PDE as an input and returns a quantity of interest. 

Additionally, we wish to mention that there exists a multitude of approaches (which are not necessarily directly RBM-or NN-related) that study the approximation of the parameter-to-solution map of PPDEs. These include methods based on sparse polynomials (see for instance \cite{CohenDeVoreHighPDE,SchwabAnalytic} and the references therein), tensors (see for instance \cite{CohenDahmenBachmayr,eigel2018variational} and the references therein) and compressed sensing (see for instance \cite{Rauhut2014CompressiveSP,DexterMixedL1} and the references therein).

Parametric PDEs also appear in the context of stochastic PDEs or PDEs with random coefficients (see for instance \cite{UQBook}) and have been theoretically examined under the perspective of \emph{uncertainty quantification}. For the sake of brevity, we only mention \cite{CohenDeVoreHighPDE} and the references therein.

Finally, we mention that a comprehensive numerical study analyzing to what extent the approximation theoretical findings of NNs (not in the context of PPDEs) are visible in practice has been carried out in \cite{adcock2020gap}. Similarly, in \cite{fokina2019growing}, a numerical algorithm that reproduces certain approximation-theoretically established exponential convergence rates of NNs was studied. The approximation rates of \cite{boelcskeiNeural} were also numerically reproduced in that paper.

\subsection{Outline}

We start by describing the parametric diffusion equation and how we discretize it in Section \ref{eq:NumericalAndTheoreticalSolutionOfPPDEs}. 
Then, we provide a formal introduction to NNs and a review of the approximation-theoretical results of NNs for parameter-to-solution maps in Section \ref{sec:ApproxByDNNs}. In Section \ref{sec:NumericalRes}, we describe our numerical experiment. We start by stating three hypotheses underlying the examples in Subsection \ref{subsec:Hypotheses}, before describing the set-up of our experiments in Subsection \ref{subsec:Setup}. After that, we present the results of the experiments in Subsection \ref{subsec:Results}. Finally, in Subsection \ref{subsec:Eval}, we evaluate and interpret the observations. In Appendix \ref{sec:EliminationCauses} we describe the measures taken to ensure comparability between test-cases.

\section{The Parametric Diffusion Equation}\label{eq:NumericalAndTheoreticalSolutionOfPPDEs}

In this section, we will introduce the abstract setup and necessary notation that we will consider throughout this paper. First of all, we will introduce the \emph{parameter-dependent diffusion equation} in Section \ref{subsec:ParDiffEq}. Afterwards, in Section \ref{subsec:HighFid}, we recapitulate some basic facts about \emph{high-fidelity discretizations} and introduce the \emph{discretized parameter-to-solution map}.

\subsection{The Parametric Diffusion Equation}\label{subsec:ParDiffEq}

Throughout this paper, we will consider the  \emph{parameter-dependent diffusion equation} with homogeneous Dirichlet boundary conditions
\begin{align}\label{eq:opeq}
    -\nabla\cdot (a(\mathbf{x}) \cdot \nabla u_a(\mathbf{x})) = f(\mathbf{x}), \quad \text{ on } \Omega= (0,1)^2, \quad u|_{\partial  \Omega} = 0,
\end{align}
where $f\in L^2(\Omega)$ is the parameter-independent right-hand side, $a\in \mathcal{A}\subset L^{\infty}(\Omega),$ and $\mathcal{A}$ constitutes some compact  set  of \emph{parametrized diffusion coefficients}. In the following, we will examine different varieties of parametrized diffusion coefficient sets $\mathcal{A}$. Following \cite{CohenDeVoreHighPDE} (by restricting ourselves to the case of finite-dimensional parameter spaces), we will always describe the elements of $\mathcal{A}$ by elements in $\R^p$ for some $p\in \N.$ 
To be more precise, we will assume that 
\begin{align} \label{eq:ParameterSetGeneral}
    \mathcal{A}= \left\{a_y:~y\in\Ycal \right\},
\end{align}
where $\Ycal\subset \R^p$ is the \emph{compact parameter space}.

A common assumption on the set $\mathcal{A}$, present in the first test-cases which we will describe below and especially convenient for the theoretical analysis of the problem, is given by \emph{affine parametrizations} of the form
\begin{align}\label{eq:Affine}
    \mathcal{A} = \left\{ a_y= a_0 + \sum_{i=1}^p y_i a_i : y=(y_i)_{i=1}^p \in \mathcal{Y}\right\},
\end{align}
where the functions $(a_i)_{i=0}^p\subset L^{\infty}(\Omega)$ are fixed.

After reparametrization, we consider the following problem, given in its variational formulation: 
\begin{align} \label{eq:paropeq}
    b_y\left(u_y, v\right) =  \int_{\Omega}f(\mathbf{x})v(\mathbf{x})~\mathrm{d}\mathbf{x}, \quad \text{ for all } y\in \Ycal,~v\in \Hil,
\end{align}
where 
\begin{align*}
    b_y:\Hil\times \Hil\to \R,~ (u,v)\mapsto\int_{\Omega} a_y(\mathbf{x}) \nabla u(\mathbf{x}) \nabla v(\mathbf{x})~\mathrm{d}\mathbf{x},
\end{align*}
and $u_y\in \Hil\coloneqq H_0^1(\Omega)$ is the solution.\footnote{Throughout this paper, we denote by $\Hil$ the space $H_0^1(\Omega)\coloneqq \{u\in H^1(\Omega):~u|_{\partial\Omega}=0\}$, where  $H^1(\Omega)\coloneqq W^{1,2}(\Omega)$ is the \emph{first-order Sobolev space} and where $\partial\Omega$ denotes the \emph{boundary of } $\Omega$. On this space, we consider the norm $\|u\|_{\Hil}=\|u\|_{H_0^1(\Omega)}\coloneqq \|u\|_{H^1(\Omega)} = \left( \sum_{|\mathbf{a}|\leq 1} \|D^{\mathbf{a}}u \|_{L^2(\Omega)}^2\right)^{1/2}.$}

We will consider experiments in which the involved bilinear forms are \emph{uniformly continuous} and \emph{uniformly coercive} in the sense that there exist $C_{\mathrm{cont}},C_{\mathrm{coer}} > 0$ with
    \begin{align*}
        \left|b_y(u,v) \right| \leq C_{\mathrm{cont}}\|u\|_{\Hil}\|v\|_{\Hil},\quad \inf_{u\in \Hil\setminus\{0\}} \frac{b_y(u,u)}{\|u\|_{\Hil}^2}  \geq C_{\mathrm{coer}}, \quad \text{ for all } u,v\in \Hil,~y\in \Ycal.
        \end{align*}
By the Lax-Milgram lemma  (see \cite[Lemma 2.1]{QuarteroniIntro}), the problem of \eqref{eq:paropeq} is \emph{well-posed}, i.e., for every $y\in \Ycal$ there exists exactly one $u_y\in \Hil $ such that \eqref{eq:paropeq} is satisfied and $u_y$ depends continuously on $f$.

\subsection{High-Fidelity Discretizations}\label{subsec:HighFid}

In practice, one cannot hope to solve \eqref{eq:paropeq} exactly for every $y\in \Ycal$. Instead, 
if we assume for the moment that $y$ is fixed,
a common approach towards the calculation of an approximate solution of \eqref{eq:paropeq} is given by the \emph{Galerkin method}, which we will describe briefly below following \cite[Appendix A]{CertReduced} and \cite[Chapter 2.4]{QuarteroniIntro}. In this framework, instead of solving \eqref{eq:paropeq}, one solves a discrete scheme of the form 
\begin{align}\label{eq:discopeq}
    b_y\left(u^{\mathrm{h}}_y , v\right) = \int_{\Omega}f(\mathbf{x})v(\mathbf{x})~\mathrm{d}\mathbf{x} \qquad \text{ for all } v\in U^{\mathrm{h}},
\end{align}
where $U^{\mathrm{h}}\subset{\Hil}$ is a subspace of ${\Hil}$ with $\mathrm{dim}\left(U^{\mathrm{h}} \right) < \infty$ and $u^{\mathrm{h}}_y\in U^{\mathrm{h}}$ is the solution of \eqref{eq:discopeq}. Let us now assume that $U^\mathrm{h}$ is given. Moreover, let $D\coloneqq \mathrm{dim}\left(U^{\mathrm{h}}\right)$, and let $\left(\varphi_{i} \right)_{i=1}^D$ be a basis for $U^{\mathrm{h}}$. Then the \emph{stiffness matrix} 
$\mathbf{B}^{\mathrm{h}}_y\coloneqq (b_y( \varphi_{j},\varphi_{i}) )_{i,j=1}^D$
is non-singular and positive definite. The solution $u^{\mathrm{h}}_y$ of \eqref{eq:discopeq} satisfies 
\begin{align*}
    u^{\mathrm{h}}_y= \sum_{i=1}^D (\mathbf{u}^{\mathrm{h}}_y)_{i} \varphi_{i},
\end{align*}
where $
    \mathbf{u}^{\mathrm{h}}_y\coloneqq (\mathbf{B}^{\mathrm{h}}_y)^{-1}\mathbf{f}^{\mathrm{h}}_y \in \R^D$ and $\mathbf{f}^{\mathrm{h}}_y\coloneqq \left( \int_{\Omega}f(\mathbf{x})\varphi_i(\mathbf{x})~\mathrm{d}\mathbf{x} \right)_{i=1}^D\in \R^D$. By Cea's Lemma (see \cite[Lemma 2.2.]{QuarteroniIntro}), $u^{\mathrm{h}}_y$ is, up to a universal constant, a best approximation of $u_y$ in $U^{\mathrm{h}}.$ 

In this framework, we can now define the central object of interest which is the map taking an element from the parameter space $\mathcal{Y}$ to the discretized solution $\mathbf{u}^{\mathrm{h}}_y$.

\begin{definition} \label{def:solutionMap}
Let $\Omega = (0,1)^2$, $U^{\mathrm{h}} \subset \Hil$ be a finite dimensional space, $\mathcal{A}\subset L^{\infty}(\Omega)$ with $\mathcal{Y} \subset \R^p$ for $p \in \N$ be as in \eqref{eq:ParameterSetGeneral}. 
Then we define the \emph{discretized parameter-to-solution map (DPtSM)} by 
\begin{align*}
\Sol \colon \mathcal{Y} \to \R^D, \qquad
y \mapsto \Sol(y) \coloneqq \mathbf{u}^{\mathrm{h}}_y.
\end{align*}
\end{definition}
\begin{remark}
The DPtSM $\Sol$ is a potentially nonlinear map from a $p$-dimensional set to a $D$-dimensional space. Therefore, without using the information that $\Sol$ has a very specific structure described through $\mathcal{A}$ and the PDE \eqref{eq:opeq}, a direct approximation of $\Sol$ as a high-dimensional smooth function will suffer from the curse of dimensionality, \cite{bellman1952theory, novak2009approximation}.
\end{remark}

Before we continue, let us introduce some crucial notation. Later, we need to compute the Sobolev norms of functions $v\in \Hil.$ This will be done via a vector representation $\mathbf{v}$ of $v$ with respect to the high-fidelity basis $(\varphi_i)_{i=1}^D.$ We denote by $\mathbf{G}\coloneqq \left(\langle \varphi_i,\varphi_j\rangle_{\Hil}\right)_{i,j=1}^{D}\in \R^{D\times D}$ the symmetric, 
positive definite \emph{Gram matrix} of the basis functions $(\varphi_i)_{i=1}^D.$ 
Then, for any $v\in U^{\mathrm{h}}$ with coefficient vector $\mathbf{v}$ with respect to the basis $(\varphi_i)_{i=1}^{D}$ we have\footnote{In this paper, $|\mathbf{x}|$ denotes the \emph{Euclidean norm} of $\mathbf{x}\in \R^n$. } (see \cite[Equation 2.41]{QuarteroniIntro})
$
\abs{\mathbf{v}}_{\mathbf{G}}\coloneqq \abs{\mathbf{G}^{1/2}\mathbf{v}}= \norm{v}_ {\Hil}.
$
In particular, $\norm{u_y^{\mathrm{h}}}_{{\Hil}} = |\mathbf{u}^{\mathrm{h}}_y |_{\mathbf{G}}$, for all $y\in \Ycal$.

\section{Approximation of the Discretized Parameter-to-Solution Map by Realizations of Neural Networks}\label{sec:ApproxByDNNs}

In this section, we describe the approximation-theoretical motivation for the numerical study performed in this paper. We present a formal definition of NNs below. In Question \ref{qu:1}, we present the underlying approximation-theoretical question of the considered learning problem. Thereafter, we recall the results of \cite{NNParametric} showing that one can upper bound the approximation rates that NNs obtain when approximating the DPtSM through an implicit notion of complexity of the DPtSM.

\subsection{Neural Networks}

NNs describe functions of compositional form that result from repeatedly applying affine linear maps and a so-called activation function. From an approximation-theoretical point of view, it is sensible to count the number of active parameters of a NN. 
To associate a meaningful and mathematically precise notion of the number of parameters to a NN, we differentiate here between \emph{neural networks} which are sets of matrices and vectors, essentially describing the parameters of the NN, and \emph{realizations of neural networks} which are the associated functions. Concretely, we make the following definition:

\begin{definition}\label{def:NeuralNetworks}
Let $n, L\in \N$.
A \emph{neural network $\Phi$ with input dimension $n$ and $L$ layers}
is a sequence of matrix-vector tuples
\[
  \Phi = \big( (\mathbf{A}_1,\mathbf{b}_1), (\mathbf{A}_2,\mathbf{b}_2), \dots, (\mathbf{A}_L, \mathbf{b}_L) \big),
\]
where $N_0 = n$ and $N_1, \dots, N_{L} \in \N$, and where each
$\mathbf{A}_\ell$ is an $N_{\ell} \times N_{\ell-1}$ matrix,
and $\mathbf{b}_\ell \in \R^{N_\ell}$.

If $\Phi$ is a NN as above, $K \subset \R^n$, and if
$\varrho\colon \R \to \R$ is arbitrary, then we define the associated
\emph{realization of $\Phi$ with activation function $\varrho$ over $K$}
(in short, the $\varrho$\emph{-realization of $\Phi$ over $K$})
as the map $\Realization_{\varrho}^K(\Phi)\colon K \to \R^{N_L}$ such that $\Realization_{\varrho}^K(\Phi)(\mathbf{x}) = \mathbf{x}_L$, where $\mathbf{x}_L$ results from the following scheme:
\begin{equation*}
  \begin{split}
    \mathbf{x}_0 &\coloneqq \mathbf{x}, \\
    \mathbf{x}_{\ell} &\coloneqq \varrho(\mathbf{A}_{\ell} \, \mathbf{x}_{\ell-1} + \mathbf{b}_\ell),
    \qquad \text{ for } \ell = 1, \dots, L-1,\\
    \mathbf{x}_L &\coloneqq \mathbf{A}_{L} \, \mathbf{x}_{L-1} + \mathbf{b}_{L},
  \end{split}
\end{equation*}
and where $\varrho$ acts componentwise, that is,
$\varrho(\mathbf{v}) \coloneqq (\varrho({v}_1), \dots, \varrho({v}_m))$
for all $\mathbf{v} = ({v}_1, \dots, {v}_s) \in \R^s$.

We call $N(\Phi) \coloneqq n + \sum_{j = 1}^L N_j$ the
\emph{number of neurons of the NN} $\Phi$ and $L$ the
\emph{number of layers}. We call $M(\Phi) \coloneqq \sum_{\ell=1}^L \|\mathbf{A}_\ell\|_0 + \|\mathbf{b}_\ell\|_0$ the \emph{number of non-zero weights of $\Phi$}.
Moreover, we refer to $N_L $ as 
\emph{the output dimension} of $\Phi$. Finally, we refer to $(N_0,\dots,N_L)$ as the \emph{architecture} of $\Phi.$
\end{definition}

We consider the following family of activation functions:

\begin{definition}\label{def:LReLU}
For $\alpha\in [0,1)$, we define by 
$\varrho_\alpha(x) \coloneqq \max\{x,\alpha x\}$
the \emph{$\alpha$-leaky rectified linear unit ($\alpha$-LReLU)}. The activation function $\varrho_0 =  \max \{x,0\}$ is called the \emph{rectified linear unit (ReLU)}.
\end{definition}
\begin{remark}\label{eq:EquivalenceOfLReLUsInApproximation}
For every $\alpha \in (0,1)$ it holds that for all $x \in \R$
$$
    \varrho_0(x) = \frac{1}{1-\alpha^2}\left(\varrho_\alpha(x) + \alpha\varrho_\alpha(-x)\right) \text { and } \varrho_\alpha(x) = \varrho_0(x) - \alpha \varrho_0(-x).
$$
Hence, for every $\alpha \in (0,1)$, we can represent the ReLU as the sum of two rescaled $\alpha$-LReLUs and vice versa. If we define for $n\in \N$
\begin{align*}
    \mathbf{P}_n(\mathbf{x}) &\coloneqq (x_1, - x_1, x_2, - x_2,\dots, x_n , - x_n) , \text{ for }  \mathbf{x}=(x_1,\dots,x_n) \in \R^{n},\\
    \mathbf{Q}_{n, \alpha}(\mathbf{x}) &\coloneqq (x_1 -\alpha  x_2, x_3 -\alpha x_4,\dots, x_{2n-1}- \alpha x_{2n}), \text{ for }  \mathbf{x}=(x_1,\dots,x_{2n}) \in \R^{2n},\\
    \mathbf{T}_{n, \alpha}(\mathbf{x}) &\coloneqq \frac{1}{1-\alpha^2}(x_1 +\alpha  x_2, x_3 +\alpha x_4,\dots, x_{2n-1} + \alpha x_{2n}), \text{ for }  \mathbf{x}=(x_1,\dots,x_{2n}) \in \R^{2n}, 
\end{align*}
then, for NNs 
\begin{align*}
\Phi_1 \! &= \! \big( (\mathbf{A}_1,\mathbf{b}_1), (\mathbf{A}_2,\mathbf{b}_2), \dots, (\mathbf{A}_L, \mathbf{b}_L) \big),\\
\Phi_2 \! &= \! \big( (\mathbf{P}_{N_1}\mathbf{A}_1, \mathbf{P}_{N_1}\mathbf{b}_1), (\mathbf{P}_{N_2}\mathbf{A}_2 \mathbf{Q}_{N_1, \alpha}, \mathbf{P}_{N_2} \mathbf{b}_2), \dots, (\mathbf{P}_{N_{L-1}}\mathbf{A}_{L-1} \mathbf{Q}_{N_{L-2}, \alpha}, \mathbf{P}_{N_{L-1}}\mathbf{b}_{L-1}), (\mathbf{A}_L \mathbf{Q}_{N_{L-1}, \alpha}, \mathbf{b}_L) \big),\\
\Phi_3 \! &= \!\big( (\mathbf{P}_{N_1}\mathbf{A}_1, \mathbf{P}_{N_1}\mathbf{b}_1), (\mathbf{P}_{N_2}\mathbf{A}_2 \mathbf{T}_{N_1, \alpha}, \mathbf{P}_{N_2} \mathbf{b}_2), \dots, (\mathbf{P}_{N_{L-1}}\mathbf{A}_{L-1} \mathbf{T}_{N_{L-2}, \alpha}, \mathbf{P}_{N_{L-1}}\mathbf{b}_{L-1}), (\mathbf{A}_L \mathbf{T}_{N_{L-1}, \alpha}, \mathbf{b}_L) \big),
\end{align*}
we have that for $K \subset \R^n$ it holds that $\Realization_{\varrho_0}^K(\Phi_1) = \Realization_{\varrho_\alpha}^K(\Phi_3)$ and $ \Realization_{\varrho_\alpha}^K(\Phi_1) = \Realization_{\varrho_0}^K(\Phi_2)$.
Moreover, it is not hard to see that $M(\Phi_1) \leq M(\Phi_2), M(\Phi_3)$ and $M(\Phi_2), M(\Phi_3) \leq 4 M(\Phi_1)$. 
Therefore, we have that for  every $\alpha_1, \alpha_2 \in [0,1)$ and every function $f: \R^n \to \R^{N_L}$ of a function space $X$ such that 
$$
    \left\| f - \Realization_{\varrho_{\alpha_1}}^K(\Phi) \right\|_{X} \leq \eps
$$
for a NN $\Phi$ implies that there exists another NN $\widetilde{\Phi}$ with $L(\widetilde{\Phi}) = L(\Phi)$ and $M(\widetilde{\Phi}) \leq 16 M(\Phi)$ such that 
$$
    \left\| f - \Realization_{\varrho_{\alpha_2}}^K(\Phi) \right\|_{X} \leq \eps.
$$
In other words, \emph{up to a multiplicative constant the parameter $\alpha$ of the $\alpha$-LReLU does not influence the approximation properties of realizations of NNs.}
\end{remark}
\begin{remark}
While Remark \ref{eq:EquivalenceOfLReLUsInApproximation} shows that all $\alpha$-LReLUs yield, in principle, the same approximation behavior, these activation functions still display quite different behavior during the training phase of NNs, where a non-vanishing parameter $\alpha$ can help avoid the occurrence of dead neurons. 
\end{remark}

\subsection{Approximation of the Discretized Parameter-to-Solution Map by Realizations of Neural Networks}\label{sec:ApproxOfSolutionMapTheoretical}

We can quantify the capability of NNs to represent the DPtSM by answering the following question:

\begin{Question}\label{qu:1}
Let $p, D \in \N$, $\alpha \in [0,1)$, $\Omega = (0,1)^2$, $U^{\mathrm{h}} \subset \Hil$ be a $D$-dimensional space, $\mathcal{A}=\{a_y:~y\in\Ycal\} \subset L^{\infty}(\Omega)$ be compact with $\mathcal{Y} \subset \R^p$ as in \eqref{eq:ParameterSetGeneral}. We consider the following equivalent questions:

\begin{itemize}
    \item For $\eps > 0$, how large do $M_\eps, L_\eps \in \N$ need to be to guarantee, that there exists a NN $\Phi$ that satisfies
\begin{itemize}
    \item[(1)] $\sup_{y \in \mathcal{Y}}|\Sol(y) - \Realization^{\mathcal{Y}}_{\varrho_\alpha}(\Phi)(y)|_{\mathbf{G}} \leq \eps,$
    \item[(2)] $M(\Phi),N(\Phi) \leq M_\eps$ and  $L(\Phi) \leq L_\eps$?
\end{itemize}
\item  For $M, L \in \N$, how small can $\eps_{L,M}>0$ be chosen so that there exists a NN $\Phi$ that satisfies
\begin{itemize}
    \item[(1)] $\sup_{y \in \mathcal{Y}}|\Sol(y) - \Realization^{\mathcal{Y}}_{\varrho_\alpha}(\Phi)(y)|_{\mathbf{G}} \leq \eps_{L,M},$
    \item[(2)] $M(\Phi),N(\Phi) \leq M$ and  $L(\Phi) \leq L$?
\end{itemize}

\end{itemize}
\end{Question}
\begin{remark}
\begin{itemize}
    \item[(i)] Conditions (1) in both instances of Question \ref{qu:1} are trivially equivalent to $$\sup_{y\in \Ycal} \left\|\sum_{i=1}^D (\mathcal{P}(y))_j\cdot\varphi_j - \left( \Realization^{\mathcal{Y}}_{\varrho_\alpha}(\Phi)(y)\right)_j \cdot\varphi_j\right\|_{\Hil}\leq \eps \text{~~or~~} \eps_{L,M}.$$
    \item[(ii)] The results to follow measure the necessary sizes of the NNs in terms of the numbers of non-zero weights $M(\Phi)$. However, from a practical point of view, we are also interested in the number of necessary neurons $N(\Phi)$. Invoking a variation of \cite[Lemma G.1.]{PetV2018OptApproxReLU} shows that similar rates to the ones below are valid for the number of neurons $N(\Phi).$
\end{itemize}
\end{remark}

If the regularity of $\Sol$ is known, then a straight-forward bound on $M_\eps$ and $L_\eps$ can be found in \cite{YAROTSKY2017103}. Indeed, if $\Sol \in C^{s}(\mathcal{Y}; \R^D)$ with $\|\Sol\|_{C^s} \leq 1$, then one can choose 
\begin{align}\label{eq:YarotskyBound}
    M_\eps \in \mathcal{O}(D \eps^{- p / s}) \text{ and } L_\eps \in \mathcal{O}(\log_2(1/\eps)), \text{ for } \eps \to 0.
\end{align}
In other situations, e.g., if $L_\eps$ is permitted to grow faster than $\log_2(1/\eps),$ one can even replace $s$ by $2s$ in \eqref{eq:YarotskyBound}, see \cite{yarotsky2018optimal, lu2020deep}.

This rate of \eqref{eq:YarotskyBound} uses the smoothness of $\Sol$ only and does not take into account the underlying structure stemming from the PDE \eqref{eq:opeq} and the choice of $\mathcal{A}$. As a result, we find this rate to be significantly suboptimal.

In \cite{NNParametric}, 
it was showed that $\Sol$ can be approximated in the sense of Question \ref{qu:1} with 
\begin{align} \label{eq:upperboundwithRBMethods}
    M_\eps &\in \mathcal{O}\left(d(\eps)D+ \left( d(\eps)^3 \log_2(d(\eps)) + pd(\eps)^2 \right) \mathrm{polylog}_2(1/\eps)\right)  \\
    L_\eps &\in \mathcal{O}\left(\mathrm{polylog}_2(1/\eps)\right), \text{ for }\eps \to 0, \nonumber
\end{align}
where $d(\eps)$ is a certain \emph{intrinsic dimension}\footnote{derived from bounds on the Kolmogorov $N$-width of $S(\Ycal)$}  of the problem, essentially reflecting the size of a \emph{reduced basis} required to sufficiently approximate $S(\Ycal).$ In many cases, especially those discussed in this manuscript, one can theoretically establish the scaling behavior of $d(\eps)$ for $\eps \to 0$. For instance, if $\mathcal{A}$ is as in \eqref{eq:Affine}, then (see \cite[Equation (3.17)] {BachCohenKolmogorov}) 
$$
    d(\eps) \in \mathcal{O}(\log_2(1/\eps)^{p}), \text{ for } \eps \to 0.
$$
Applied to \eqref{eq:upperboundwithRBMethods} this yields that
\begin{align*}
    M_\eps \in \mathcal{O}\left(D\log_2(1/\eps)^{ p}+ p \cdot \log_2(1/\eps)^{c p}\right), \text{ for } \eps \to 0,
\end{align*}
for some $c\geq 1.$ We also mention a similar approximation result, not of the discretized parametric map $\Sol$ but of the parametrized solution $(y,\mathbf{x}) \mapsto u_{a_y}(\mathbf{x})$,
where $u_{a_y}$ is as in \eqref{eq:opeq} for $a=a_y$.
In this situation, and for specific parametrizations of $\mathcal{A}$, \cite[Theorem 4.8]{SchwabZech} shows that this map can be approximated by the realization of a NN using the ReLU activation function up to an error of $\eps$ with a number of weights that essentially scales like $\eps^{-r}$ where $r$ depends on the summability of the (in this case potentially infinite) sequence $(a_i)_{i=1}^{\infty}$ such that $a_y= a_0 + \sum_{i=1}^\infty y_i a_i$ for a coefficient vector $ y=(y_i)_{i=1}^{\infty} $. Here $r$ can be very small if $\|a_i\|_{L^\infty(\Omega)}$ decays quickly for $i \to \infty$. This leads to very efficient approximations. 

While the aforementioned results all examine the approximation-theoretical properties of realizations of NNs with respect to the \emph{uniform} approximation error, they trivially imply the same rates if we examine the \emph{average} errors 
$$\left(\int_{\Ycal} \left|\Sol(y) - \Realization^{\mathcal{Y}}_{\varrho_\alpha}(\Phi_\eps)(y) \right|_{\mathbf{G}}^p ~\mathrm{d}\mu(y)\right)^{1/p},$$ which are often used in practice. 
Here,   $1\leq p<\infty$ and $\mu$ is an arbitrary probability measure on $\Ycal.$ In this paper, we examine the discrete counterpart of the \emph{mean relative error}
$$
\int_{\Ycal} \frac{ \left|\Sol(y) - \Realization^{\mathcal{Y}}_{\varrho_\alpha}(\Phi_\eps)(y) \right|_{\mathbf{G}}}{|\Sol(y)|_{\mathbf{G}}} ~\mathrm{d}\mu(y),
$$
where $\mu$ denotes the \emph{uniform probability measure} on $\Ycal$.

In view of the aforementioned theoretical results, it is clear that a parameter that is not the dimension of the parameter space $\mathcal{Y}$ but a problem-specific notion of complexity determines the hardness of the approximation problem of Question \ref{qu:1}. To what extent this theoretical observation influences the hardness of the practical learning problem will be analyzed in the numerical experiment presented in the next section.

\section{Numerical Survey of Approximability of Discretized Parameter-to-Solution Maps} \label{sec:NumericalRes}

As outlined in Section \ref{sec:ApproxByDNNs}, the theoretical hardness of the approximation problem of Question \ref{qu:1} is determined by an intrinsic notion of complexity that potentially differs substantially from the dimension of the parameter space. 

To test how this intrinsic complexity affects the practical machine-learning based solution of \eqref{eq:opeq}, we perform a comprehensive study where we train NNs to approximate the DPtSM $\Sol$ for various choices of $\mathcal{A}$. Here, we are especially interested in the performance of the learned approximation of $\Sol$ for varying complexities of $\mathcal{A}$.
In this context, we test the hypotheses listed in the following Subsection \ref{subsec:Hypotheses}. The remainder of this section is structured as follows: In Subsection \ref{subsec:Setup}, we introduce the concrete setup of parametrized diffusion coefficient sets, NN architecture, and optimization procedure and explain how the choice of test-cases are related to our hypotheses. Afterwards, in Subsection \ref{subsec:Results}, we report the results of our numerical experiments. Subsection \ref{subsec:Eval} is devoted to an evaluation and interpretation of these results in view of the hypotheses of Subsection \ref{subsec:Hypotheses}.

\subsection{Hypotheses}\label{subsec:Hypotheses}

\begin{itemize}
    \item[\textbf{[H1]}] \textit{The performance of learning the DPtSM does not suffer from the curse of dimensionality:}

    The theoretical results of \cite{NNParametric} show that the dimension of the parameter space $p$ is not the main factor in determining the hardness of the underlying approximation-theoretical problem. As already outlined in the introduction, it is by no means clear that this effect is visible in a practical learning problem. 
    
    We expect that after accounting for effects stemming from optimization and sampling to promote comparability between test-cases in a way described in Appendix \ref{sec:EliminationCauses}, the performance of the learning method will scale only mildly with the dimension of the parameter space.
    
    \item[\textbf{[H2]}] \textit{The performance of learning the DPtSM is very sensitive to parametrization:}
    
    We expect that, within the framework of Question \ref{qu:1}, there are still extreme differences of intrinsic complexities for different choices of  parametrizations for the diffusion coefficient sets $\mathcal{A} \subset L^\infty(\Omega)$ as defined in  \eqref{eq:ParameterSetGeneral}.     
    However, it is not clear to what extent NNs are capable of resolving the low-dimensional sub-structures generated by various choices of $\mathcal{A} \subset L^\infty(\Omega)$.
    
    Since realizations of NNs are a very versatile function class, we expect the degree to which the performance of a trained NN depends on the number of parameters to vary strongly over the choice of $(a_i)_{i=1}^p$.

    \item[\textbf{[H3]}] \textit{Learning the DPtSM is efficient also for non-affinely parametrized problems:}
    
    The analysis of PPDEs often relies on affine parametrizations  as in \eqref{eq:Affine} or smooth variations thereof. 
    
    We expect the overall theme that NNs perform according to an intrinsic complexity of the problem depending only weakly on the parameter dimension to hold in more general cases.

\end{itemize}

\subsection{Setup of Experiments}\label{subsec:Setup}

To test the hypotheses \textbf{[H1]}, \textbf{[H2]}, and \textbf{[H3]} of Section \ref{subsec:Hypotheses}, we consider the following setup.

\subsubsection{Parameterized Diffusion Coefficient Sets}

We perform training of NNs for different instances of the approximation problem of Question \ref{qu:1}. Here, we always assume the right-hand side to be fixed as $f(\mathbf{x})= 20+10x_1-5x_2$, for $\mathbf{x} = (x_1,x_2) \in \Omega$, and we vary the parametrized diffusion coefficient set $\mathcal{A}$.

We consider four different parametrized diffusion coefficient sets as described in the test-cases \textbf{[T1]}-\textbf{[T4]} (for a visualization of \textbf{[T3]} and \textbf{[T4]} see Figure \ref{fig:Setups} below). \textbf{[T1]}, \textbf{[T2]} and \textbf{[T3-F]} are affinely parametrized whereas the remaining parametrizations are non-affine. 

\begin{itemize}

    \item[\textbf{[T1] Trigonometric Polynomials:}] In this case, the set $\mathcal{A}$ consists of trigonometric polynomials that are weighted according to a scaling coefficient $\sigma$. To be more precise, we consider
    \begin{align*}
        \mathcal{A}^{\mathrm{tp}}(p, \sigma) \coloneqq \left\{\mu+ \sum_{i=1}^p  y_i \cdot i^{\sigma} \cdot (1+a_i):~y\in \Ycal=[0,1]^p \right\},
    \end{align*}
    for some fixed shift $\mu>0$ and a scaling coefficient $\sigma\in \R$. Here $a_i(\mathbf{x})=\sin\left(\left\lfloor \frac{i+2}{2}\right\rfloor \pi x_1 \right) \sin\left(\left\lceil \frac{i+2}{2}\right\rceil  \pi x_2\right),$ for $i=1,\dots, p.$
    
    We analyze the cases $p=2, 5, 10, 15, 20$ and, for each $p$, the scaling coefficients $\sigma=-1,0,1$. As a shift we always choose $\mu = 1$.
        
    \item[\textbf{[T2] Chessboard Partition:}] Here, we assume that $p=s^2$ for some $s\in \N$ and we consider\footnote{$\mathcal{X}_A$ denotes the \emph{indicator function} of $A$.}
    \begin{align*}
        \mathcal{A}^{\mathrm{cb}}(p,\mu) \coloneqq \left\{\mu + \sum_{i=1}^p y_i\mathcal{X}_{\Omega_i}:~y\in \Ycal= [0,1]^p \right\}, 
    \end{align*} 
    where $(\Omega_i)_{i=1}^p$ forms a $s\times s$ chessboard partition of $(0,1)^2$ and $\mu>0$ is a fixed shift. 
      
    We examine this test-case for the shifts $\mu=10^{-1}, 10^{-2}, 10^{-3}$, and, for each $\mu$ we consider  $s=2,3,4,5$ which yields $p=4,9,16,25,$ respectively.

     \item[\textbf{[T3] Cookies:}]  In this test-case we differentiate between two sub-cases: 
     \begin{itemize} 
     \item[\textbf{[T3-F] Cookies with Fixed Radii:}]
     In this setting, we assume that $p=s^2$ for some $s\in \N$ and we consider 
    \begin{align*}
        \mathcal{A}^{\mathrm{cfr}}(p,\mu) \coloneqq \left\{\mu + \sum_{i=1}^{p} y_i\mathcal{X}_{\Omega_{i}}:~y\in \Ycal= [0,1]^p \right\},
    \end{align*}   
    for some fixed shift $\mu>0$ 
    where the  $\Omega_{i}$ are disks with centers $((2k+1)/(2s),(2\ell-1)/(2s))$, where $i= ks+\ell$ for uniquely determined $k\in\{0,\dots s-1\}$ and $\ell\in \{ 1,\dots, s\}$. The radius is set to $r/(2s)$ for some fixed $r\in (0,1]$.
    
     We examine this test-case for fixed $\mu=10^{-4},$ $r=0.8$ and $s=2,3,4,5,6$ which yields parameter dimensions $p=4,9,16,25,36,$ respectively.
     \item[\textbf{[T3-V] Cookies with Variable Radii:}]
     Here, we additionally assume that the radii of the involved disks are not fixed anymore. To be more precise, for $s\in\N$ and every $i=1,\dots,s,$ we are given disks $\Omega_{i, y_{i+s^2}}$ with center as before and radius $y_{i+s^2}/(2s)$ for $y_{i+s^2}\in  [0.5,0.9],$
      so that $\Ycal=[0,1]^{s^2}\times [0.5,0.9]^{s^2} \subset \R^{p}$ with $p=2s^2$. We define 
    \begin{align*}
        \mathcal{A}^{\mathrm{cvr}}(p,\mu) \coloneqq \left\{\mu + \sum_{i=1}^{p} y_i\mathcal{X}_{\Omega_{i, y_{i+s^2}}}:~y\in \Ycal = [0,1]^{p}\times [0.5,0.9]^{p}\right\}.
    \end{align*}
      Note that, $\mathcal{A}^{\mathrm{cvr}}(p,\mu)$ is \emph{not} an affine parametrization. 
      
      We consider the shifts $\mu=10^{-4}$ and $\mu=10^{-1},$ and, for each $\mu$, we consider the cases $s=2,3,4,5$ which yields the parameter dimensions $p= 8,18,32,50,$ respectively. 
      \end{itemize}
      
    \item[\textbf{[T4] Clipped Polynomials:}] Let
    \begin{align*}
        \mathcal{A}^{\mathrm{cp}}(p,\mu) \coloneqq \left\{ \max\left\{\mu,\sum_{i=1}^p y_i m_i\right\}:~ (y_i)_{i=1}^p \in \mathcal{Y}= [-1,1]^p \right\},
    \end{align*}
    where $\mu >0$ is the fixed clipping value and
    $(m_i)_{i=1}^p$ is the monomial basis of the space of all two-variate polynomials of degree $\leq k$. Therefore $p= \binom{2+k}{2}.$
    
    We examine this test-case for fixed shift $\mu=10^{-1}$ and for $k=2,3,5,8,12$ which yields parameter dimensions $p=6,10,21,45,91,$ respectively. 
      
\end{itemize}

\begin{figure}[H]
    \centering
    \includegraphics[width=.9\textwidth]{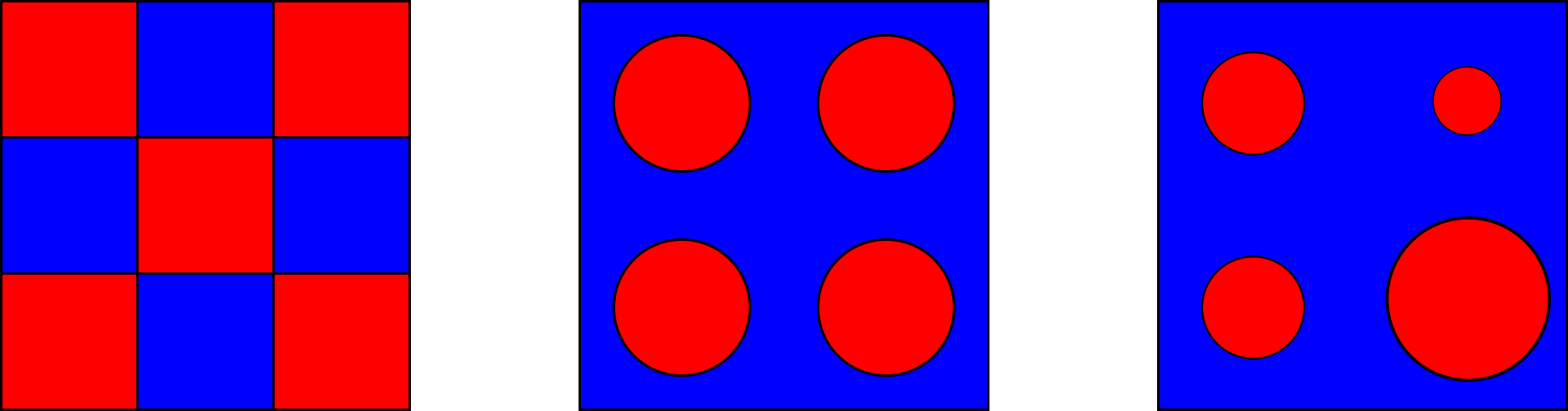}
    \caption{Partition of $\Omega$ as in  Test-case \textbf{[T2]} (left) for $p=9,$ (the red and blue areas indicate the $\Omega_i$), test-case \textbf{[T3-F]} (middle) for $p=4$ (the red areas indicate the $\Omega_{i}$) and test-case \textbf{[T3-V]} (right) for $p=8$ (the red areas indicate the $\Omega_{i,y_{i+s^2}}$).} 
    \label{fig:Setups}
\end{figure}

\subsubsection{Setup of Neural Networks and Training Procedure}\label{sec:SetupNN}

Our experiments are implemented using Tensorflow, \cite{Tensorflow}, for the learning procedure and FEniCS, \cite{Fenics}, as FEM solver. The code used for dataset generation of all considered test-cases is made publicly available at \url{www.github.com/MoGeist/diffusion_PPDE}. To be able to compare different test-cases and remove all effects stemming from the optimization procedure, we train almost the same model for all parameter spaces. The only---to a certain extent inevitable---change that we allow between test-cases is that the input dimension of the NN changes to that of the parameter space. Concretely, we consider the following setup:

\begin{itemize}
    \item[(1)] The \emph{finite element space} $U^{\mathrm{h}}$ resulting from a triangulation of $\Omega=[0,1]^2$ with $101\times 101= 10201$ equidistant grid points and first-order Lagrange finite elements. This space shall serve as a discretized version of the space $H^1(\Omega).$  We denote by $D=10201$ its dimension and by $(\varphi_i)_{i=1}^D$ the corresponding finite element basis.
    \item[(2)] The (feedforward) \emph{neural network architecture} $S=(p,300,\dots,300,10201)$ with $L=11$ layers, where $p$ is test-case-dependent and the weights and biases are initialized according to a normal distribution with mean $0$ and standard deviation $0.1$. 
    \item[(3)] The \emph{activation function} is the 0.2-LReLU of Definition \ref{def:LReLU}.
    \item[(4)] The \emph{loss function} is the relative error on the finite-element discretization of $\Hil$
    $$
        \mathcal{L}:\R^D\times (\R^D\setminus\{0\})\to \R,\quad (\mathbf{x}_1,\mathbf{x}_2)\mapsto \frac{\abs{\mathbf{x}_1-\mathbf{x}_2}_{\mathbf{G}}}{\abs{\mathbf{x}_2}_{\mathbf{G}}}.
    $$
    \item[(5)] The  \emph{training set} $(y^{i,\mathrm{tr}})_{i=1}^{N_{\mathrm{train}}}\subset \Ycal$ consists of $N_{\mathrm{train}} \coloneqq 20000$ i.i.d. parameter samples, drawn with respect to the \emph{uniform probability measure} on $\Ycal.$ 
    \item[(6)] The \emph{test set} $(y^{i,\mathrm{ts}})_{i=1}^{N_{\mathrm{test}}}\subset \Ycal$ consists of $N_{\mathrm{test}} \coloneqq 5000$ i.i.d. parameter samples,  drawn with respect to the \emph{uniform probability measure} on $\Ycal.$ 
\end{itemize}

In our experiments, we aim at finding a NN $\Phi$ with architecture $S$ such that the \emph{mean relative training error}
\begin{align*}
    \frac{1}{N_{\mathrm{train}}} \sum_{i=1}^{N_{\mathrm{train}}}  \mathcal{L}\left(\Realization_\varrho^{\Ycal}\left(\Phi \right)(y^{i,\mathrm{tr}}),\mathbf{u}_{y^{i,\mathrm{tr}}}^{\mathrm{h}}\right) = \frac{1}{N_{\mathrm{train}}} \sum_{i=1}^{N_{\mathrm{train}}} \frac{\left\| \sum_{j=1}^{D} \left(\Realization_\varrho^{\Ycal}\left(\Phi \right)(y^{i,\mathrm{tr}}) \right)_j\cdot \varphi_j -u_{y^{i,\mathrm{tr}}}^{\mathrm{h}}\right\|_{\Hil}}{\left\|u_{y^{i,\mathrm{tr}}}^{\mathrm{h}}\right\|_{\Hil}}
\end{align*}
is minimized. We then test the accuracy of our NN 
by computing the \emph{mean relative test error}
\begin{align*}
     \frac{1}{N_{\mathrm{test}}} \sum_{i=1}^{N_{\mathrm{test}}} \mathcal{L}\left(\Realization_\varrho^{\Ycal}\left(\Phi \right)(y^{i,\mathrm{ts}}),\mathbf{u}_{y^{i,\mathrm{ts}}}^{\mathrm{h}}\right) = \frac{1}{N_{\mathrm{test}}} \sum_{i=1}^{N_{\mathrm{test}}} \frac{\left\| \sum_{j=1}^{D} \left(\Realization_\varrho^{\Ycal}\left(\Phi \right)(y^{i,\mathrm{ts}}) \right)_j\cdot \varphi_j -u_{y^{i,\mathrm{ts}}}^{\mathrm{h}}\right\|_{\Hil}}{\left\|u_{y^{i,\mathrm{ts}}}^{\mathrm{h}}\right\|_{\Hil}}.
\end{align*}
Here, we use the mean \emph{relative} error instead of the mean absolute error in order to establish comparability of our results between different sets $\mathcal{A}$, allowing us to put our results into context.

The optimization is done through \emph{batch gradient descent}. 
To ensure further comparability between the different setups, the hyper-parameters in the optimization procedure are kept fixed: Training is conducted with batches of size 256 using the ADAM optimizer \cite{ADAM} with hyper-parameters $\alpha = 2.0\times 10^{-4}$, $ \beta_{1}=0.9$, $\beta_{2}=0.999$ and $\eps=1.0\times 10^{-8}$. Training is stopped after reaching $40000$ epochs. Having trained the NN, 
for some new input $y\in\mathcal{Y},$ the computation of the approximate discretized solution $\Realization_\varrho^{\Ycal}(\Phi)(y)$ is done by a simple forward pass.

\subsubsection{Relation to Hypotheses}

\begin{figure}[htb]
\centering
  \begin{tabular}{@{}cc@{}}
    \includegraphics[width=.45\textwidth]{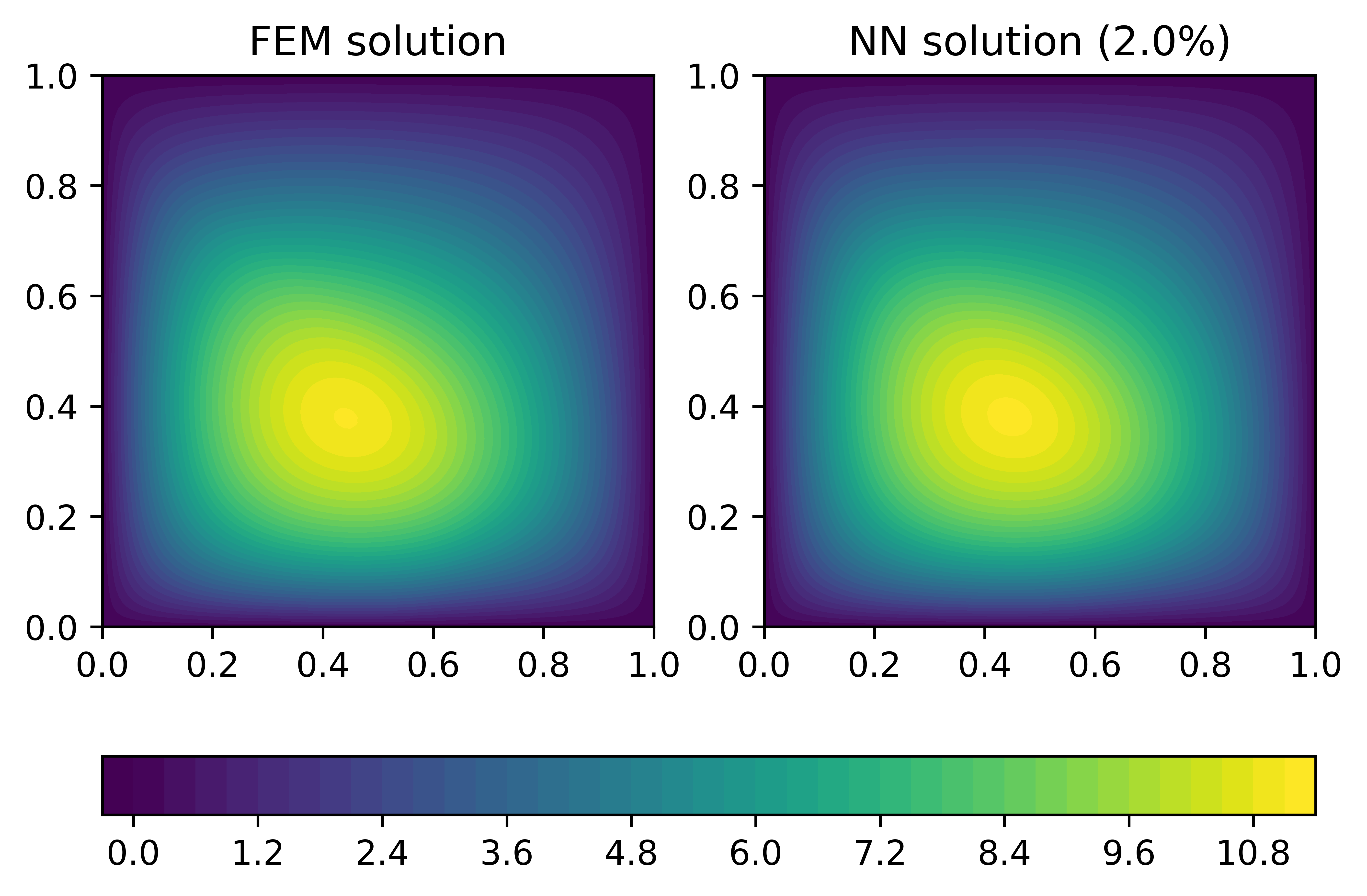} &
    \includegraphics[width=.45\textwidth]{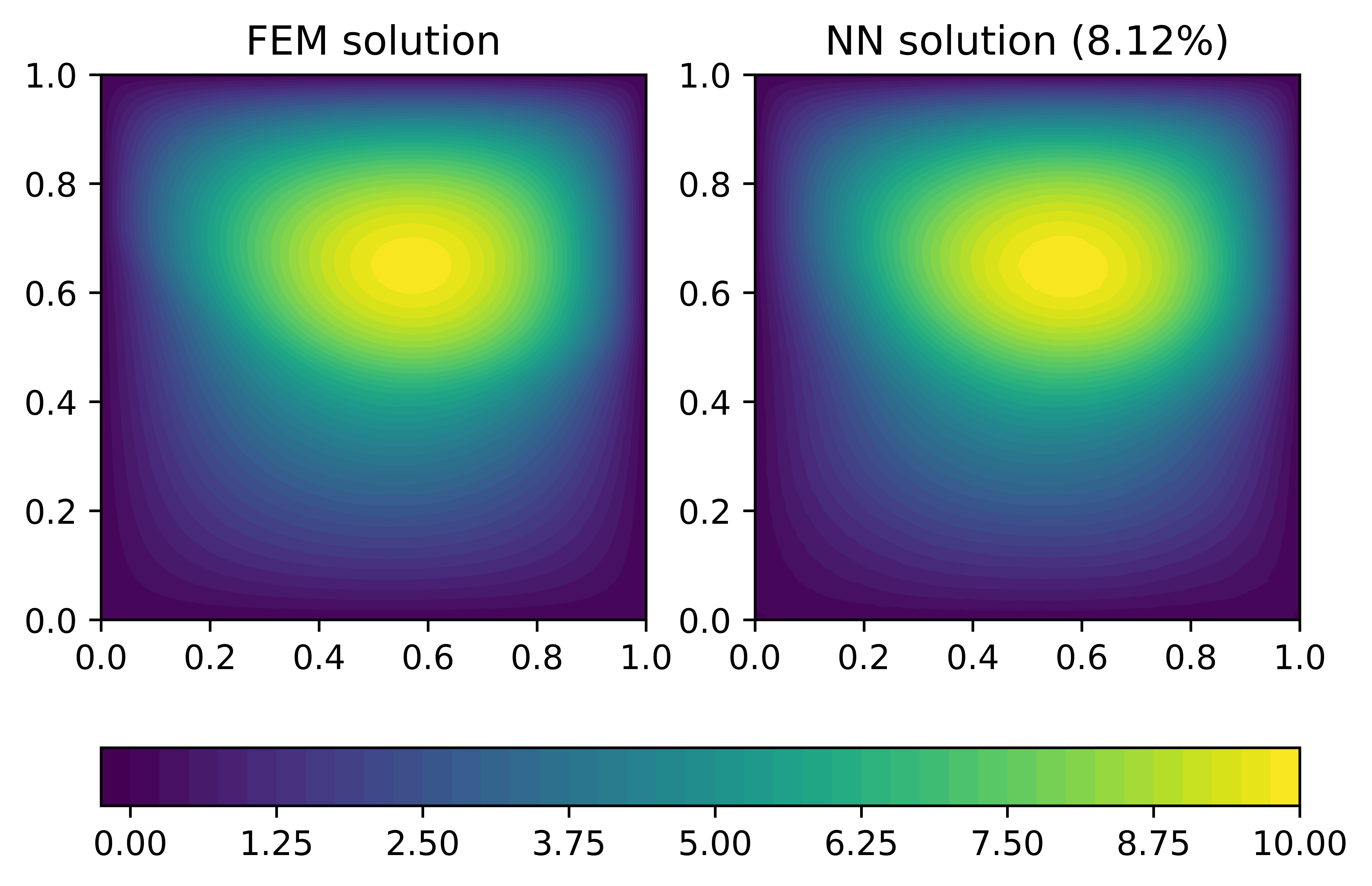} \\
    \includegraphics[width=.45\textwidth]{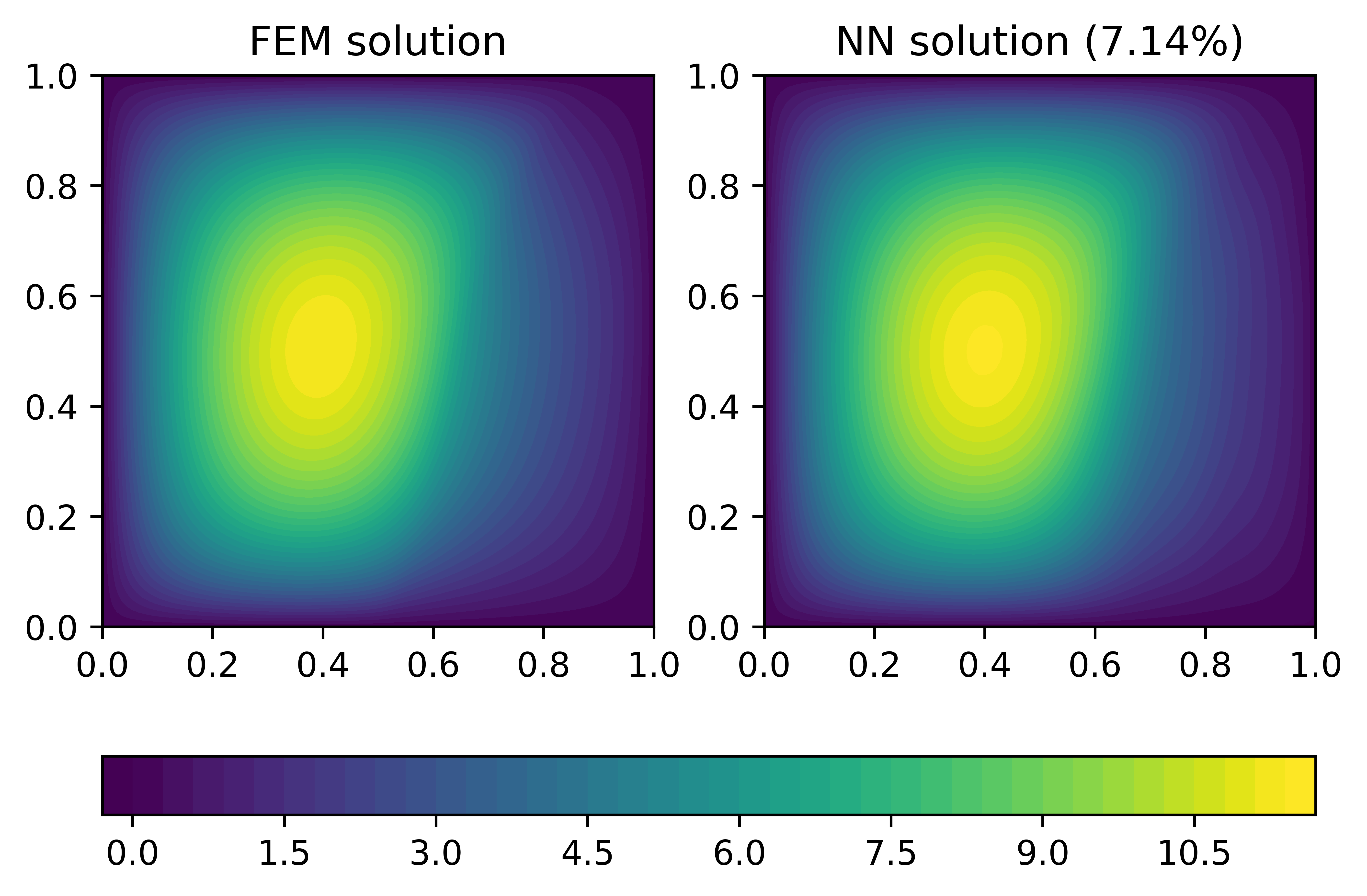} &
    \includegraphics[width=.45\textwidth]{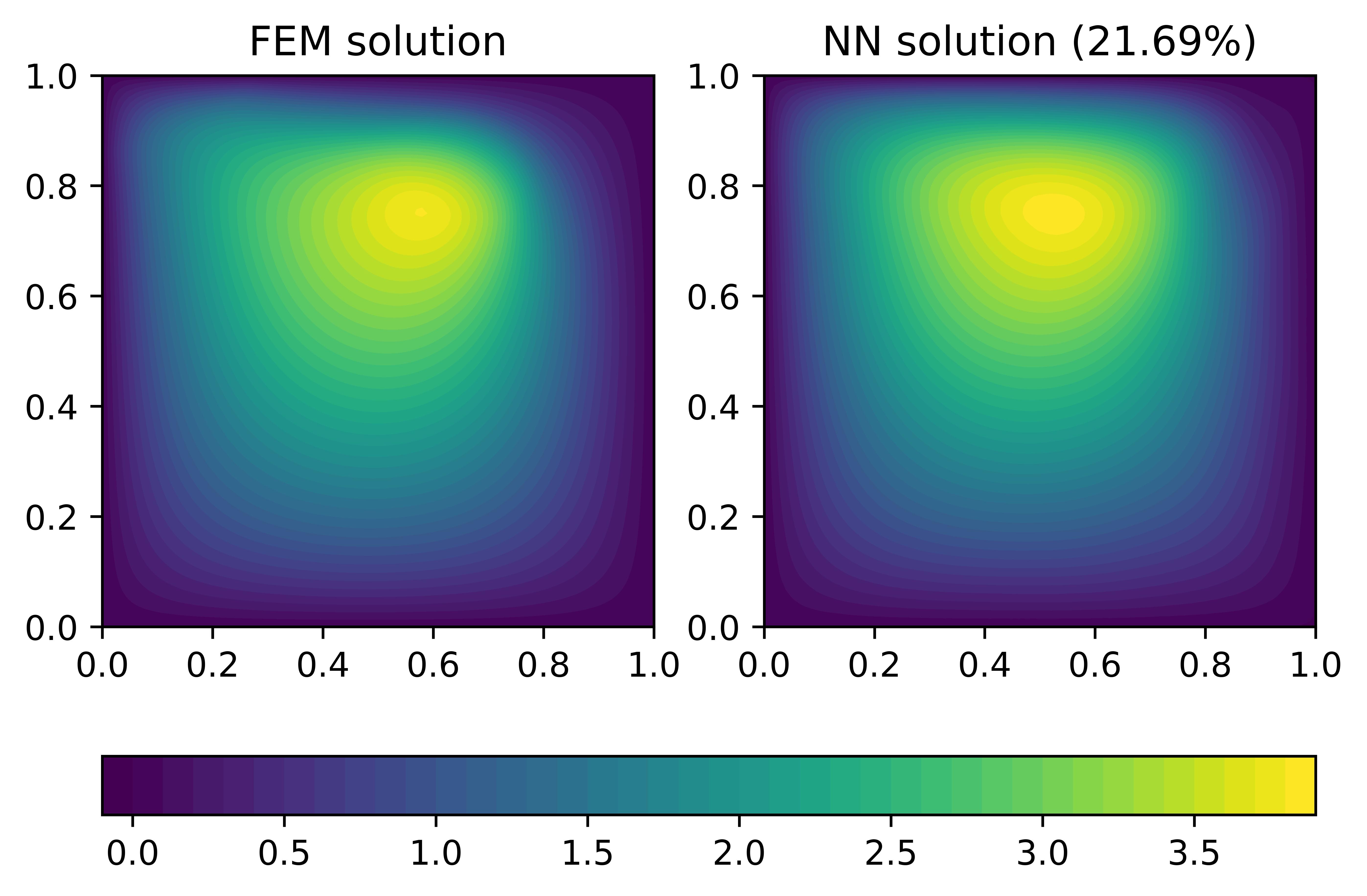}    
  \end{tabular}
  \caption{Comparison of the ground truth solution and the one generated by the NN for an average (left) and a poor performing case (right) for $\mu=10^{-1}$ and $p=6$ (top) and $p=91$ (bottom) for test-case \textbf{[T4]}. The percentage in brackets represents the relative test error for this particular sample.}\label{fig:polynomialsTest}
\end{figure}

The test-cases \textbf{[T1]} - \textbf{[T4]} are designed to test the hypotheses \textbf{[H1]} - \textbf{[H3]} in the following way:

\begin{itemize}
    \item[ \textit{Enabling comparability between test-cases:}]
    We implement three measures to produce a uniform influence of the optimization and sampling procedure in all test-cases. These are that we only change the architecture in the minimally required way between test-cases to not alter the optimization behavior, we analyze a posteriori the optimization behavior to see if there are qualitative differences between test-cases, and we choose the number of training samples in such a way that neither moderate further increasing or decreasing of the number of training samples affects the outcome of the experiments. We describe these measures in detail in Appendix \ref{sec:EliminationCauses}. 

    \item[\textit{Relation to Hypothesis} \textbf{[H1]}:] 
    
    To test if the learning method suffers from the curse of dimensionality or if the prediction of \cite{NNParametric} that its complexity is determined only by some intrinsic complexity of the function class holds, we run all test-cases \textbf{[T1]}-\textbf{[T4]} for various values of the dimension of the parameter space, and study the resulting scaling behavior.

    \item[\textit{Relation to Hypothesis} \textbf{[H2]}:] To understand the extent to which the NN model is sufficiently versatile to adapt to various types of solution sets, we study four commonly considered parametrized diffusion coefficient sets which also include multiple subproblems described via the hyper-parameters $\sigma$ and $\mu$. The parametrized sets exhibit the following different characteristics:
    
    \begin{itemize}
        \item[\textbf{[T1]}] The parameter-dependence in this case is affine (i.e. the forward-map $y\mapsto b_y(u,v)$ depends affinely on $y$ for all $u,v\in \Hil$) whereas the spatial regularity of the functions $(a_i)_{i=1}^p$ is analytic. To vary the difficulty of the problem at hand, we consider different instances of the scaling coefficient $\sigma$ which put different emphasis on the high-frequency components of the functions $(a_i)_{i=1}^p$. In particular, if $\sigma >0,$ a higher weight is put on the high-frequency components than on the low-frequency ones  whereas the opposite is true for $\sigma<0.$ 
        
        \item[\textbf{[T2]}] The parameter-dependence in this case is affine again, whereas the spatial regularity of the $(\mathcal{X}_{\Omega_i})_{i=1}^p$ is very low. To vary the difficulty of the problem, we consider different instances of shifts $\mu$. The higher the shift is, the more elliptic the problem becomes. 
        
        \item[\textbf{[T3]}] \textbf{[T3-F]} again exhibits affine parameter-dependence and the same regularity properties as test-case \textbf{[T2]}. However, this problem is considered to be easier than test-case \textbf{[T2]} since the $\overline{\Omega_i}$ do not intersect each other. 
        
        For test-case \textbf{[T3-V]}, the geometric properties of the domain partition are additionally encoded via a parameter thereby rendering the problem to be non-affine.   
        
        \item[\textbf{[T4]}] In this case, the parameter-dependence is non-affine and has low regularity due to the clipping procedure. Additionally, the spatial regularity of the functions $a_y$ is comparatively low in general.
    \end{itemize}
    
    A visualization highlighting the versatility of our test-cases can be seen when comparing the FE solutions in Figure \ref{fig:SquaresExperiment} (test-case \textbf{[T2]})   with the FE solutions in Figure \ref{fig:polynomialsTest} (test-case \textbf{[T4]}).
    
   \item[\textit{Relation to Hypothesis} \textbf{[H3]}:] The test-cases \textbf{[T3-V]} and \textbf{[T4]} are non-affinely parametrized.   
    
\end{itemize}

\subsection{Numerical Results}\label{subsec:Results}

In this subsection, we report the results of the test-cases announced in the previous subsection. 

\subsubsection*{\textbf{[T1]} Trigonometric Polynomials} 

We observe the following mean relative test errors for the sets $\mathcal{A}^{\mathrm{tp}}(p, \sigma)$.
\begin{table}[H]
\renewcommand{\arraystretch}{1.5}
\vspace{-0.5cm}
\small
\begin{tabularx}{\textwidth}{|l|l|l|l|l|l|}
  \hline
    \textbf{Parameter dimension $p$} & 2 & 5 & 10 & 15 & 20 \\ \hline
 \textbf{Mean relative test error } ($\sigma$ = -1) & 0.32 \% & 0.36 \% & 0.42 \% &  0.43 \%  & 0.43 \%
 \\ \hline
 \textbf{Mean relative test error } ($\sigma$ = 0) &  0.36 \% &  0.43 \% &  0.44 \% &   0.51 \%&   0.59 \%
 \\ \hline
  \textbf{Mean relative test error } ($\sigma$ = 1) &  0.39 \% &  0.84 \%  & 2.05 \% &   2.45 \% &  3.85 \%
 \\ \hline
  \end{tabularx}
   \caption{Mean relative test error for test-case \textbf{[T1]} and different parameter dimensions $p$, different scaling parameters $\sigma$ and shift $\mu=1$.}
   \label{tab:TrigonPol}
\end{table}

\subsubsection*{\textbf{[T2]} Chessboard Partition} 

We observe the following mean relative test errors for the sets $\mathcal{A}^{\mathrm{cb}}(p,\mu)$.
\begin{table}[H]
\renewcommand{\arraystretch}{1.5}
\vspace{-0.5cm}
\small
\begin{tabularx}{\textwidth}{|l|l|l|l|l|}
  \hline
  $s$ & 2 & 3 & 4 & 5  \\ \hline
    \textbf{Parameter dimension $p $} & 4 & 9  & 16 & 25 \\ \hline
  \textbf{Mean relative test error }($\mu = 10^{-1}$)   & 0.57 \% & 1.06 \% & 2.19 \% & 3.22 \%  \\ \hline
     \textbf{Mean relative test error} ($\mu = 10^{-2}$)  & 0.66 \% & 1.81 \% & 4.13 \% & 6.78 \%  \\ \hline
     \textbf{Mean relative test error} ($\mu = 10^{-3}$) & 1.09 \% & 4.47 \% & 12.01 \% & 23.96 \%  \\ \hline
  \end{tabularx}
   \caption{Mean relative test error for test-case \textbf{[T2]} and parameter dimensions $p=s^2$.}
     \label{tab:SquareI}
\end{table}

\begin{figure}[htb]
\centering
 \begin{tabular}{@{}cc@{}}
    \includegraphics[width=.45\textwidth]{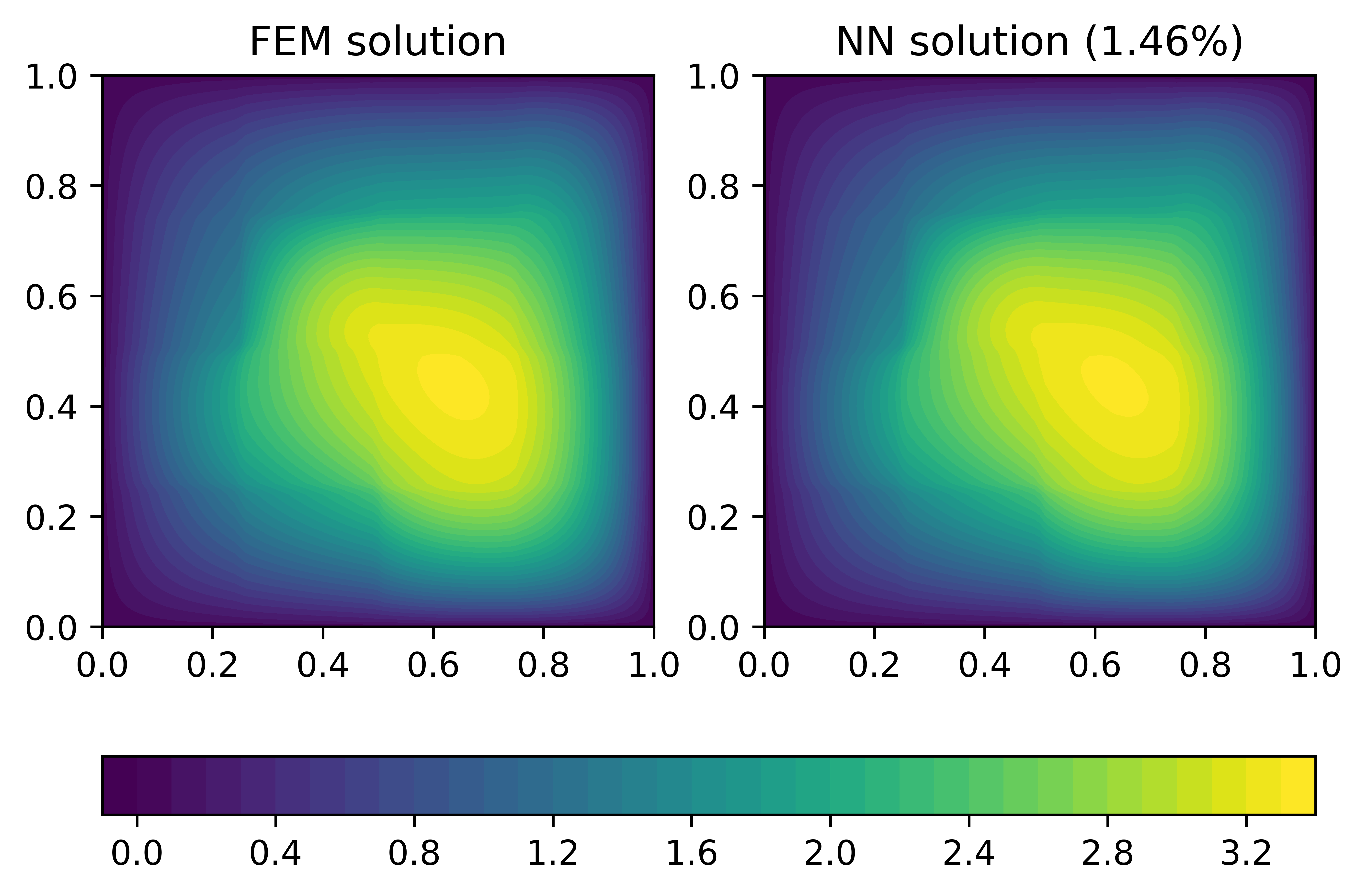} &
    \includegraphics[width=.45\textwidth]{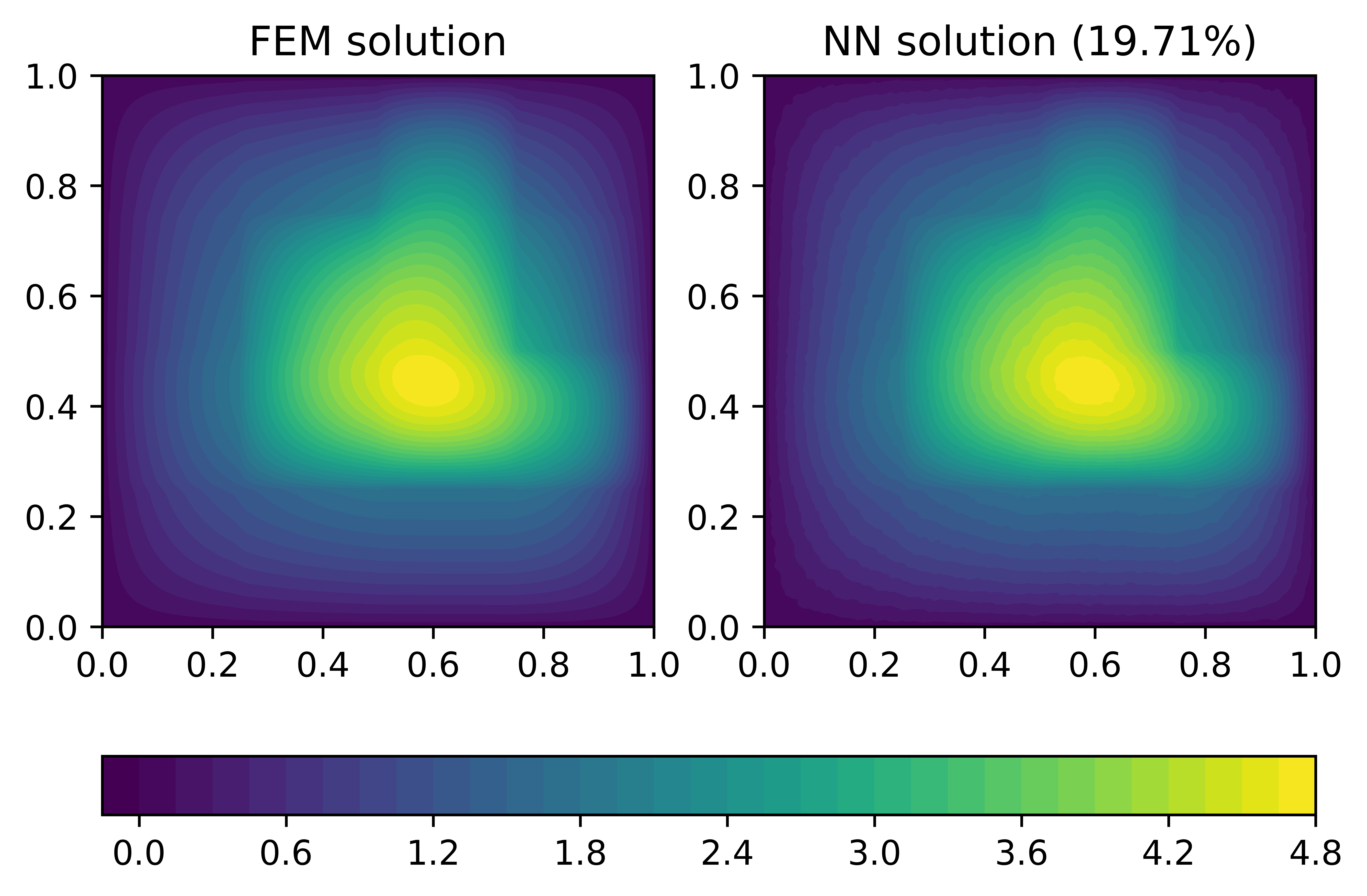} \\
    \includegraphics[width=.45\textwidth]{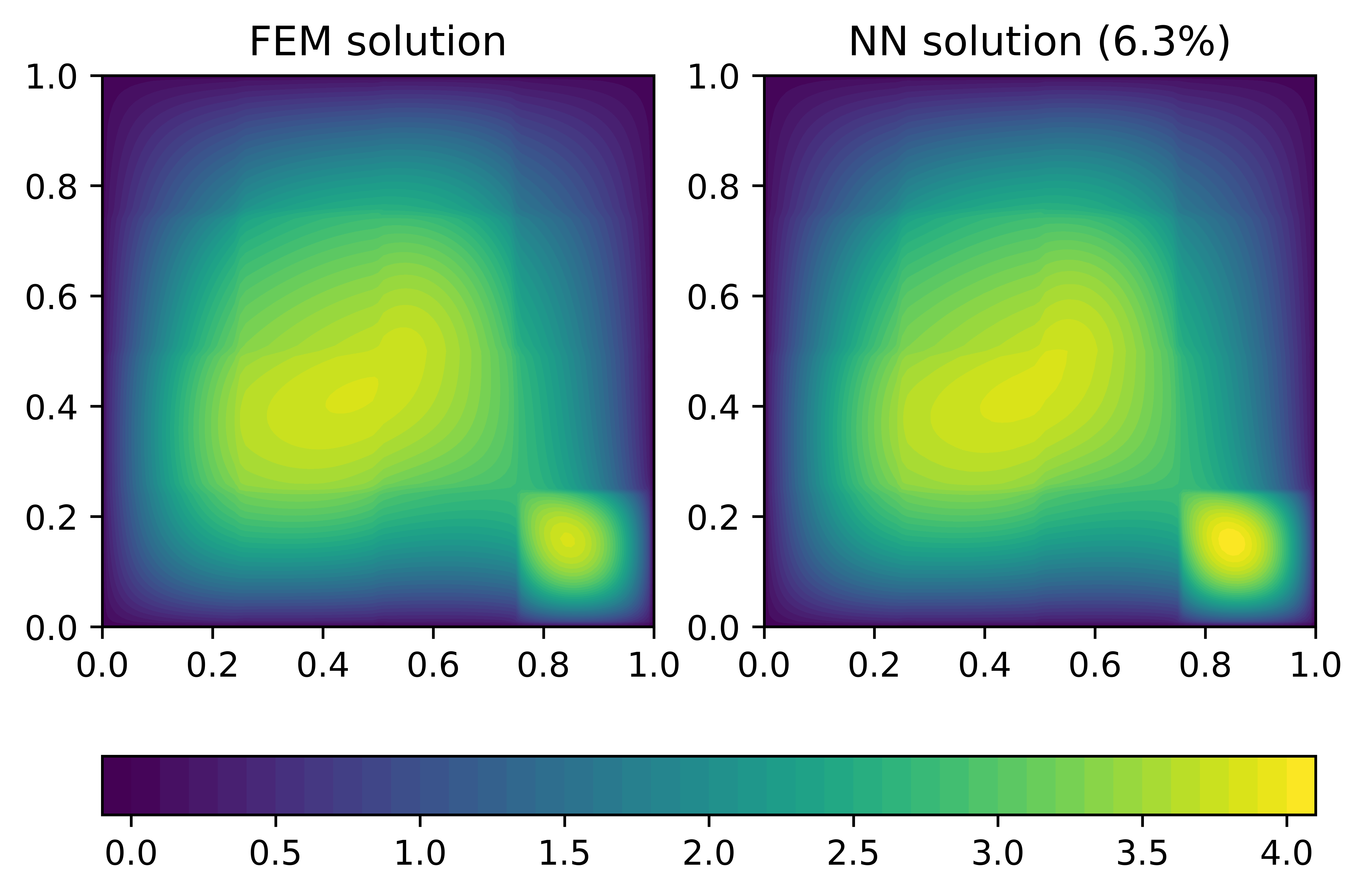} &
    \includegraphics[width=.45\textwidth]{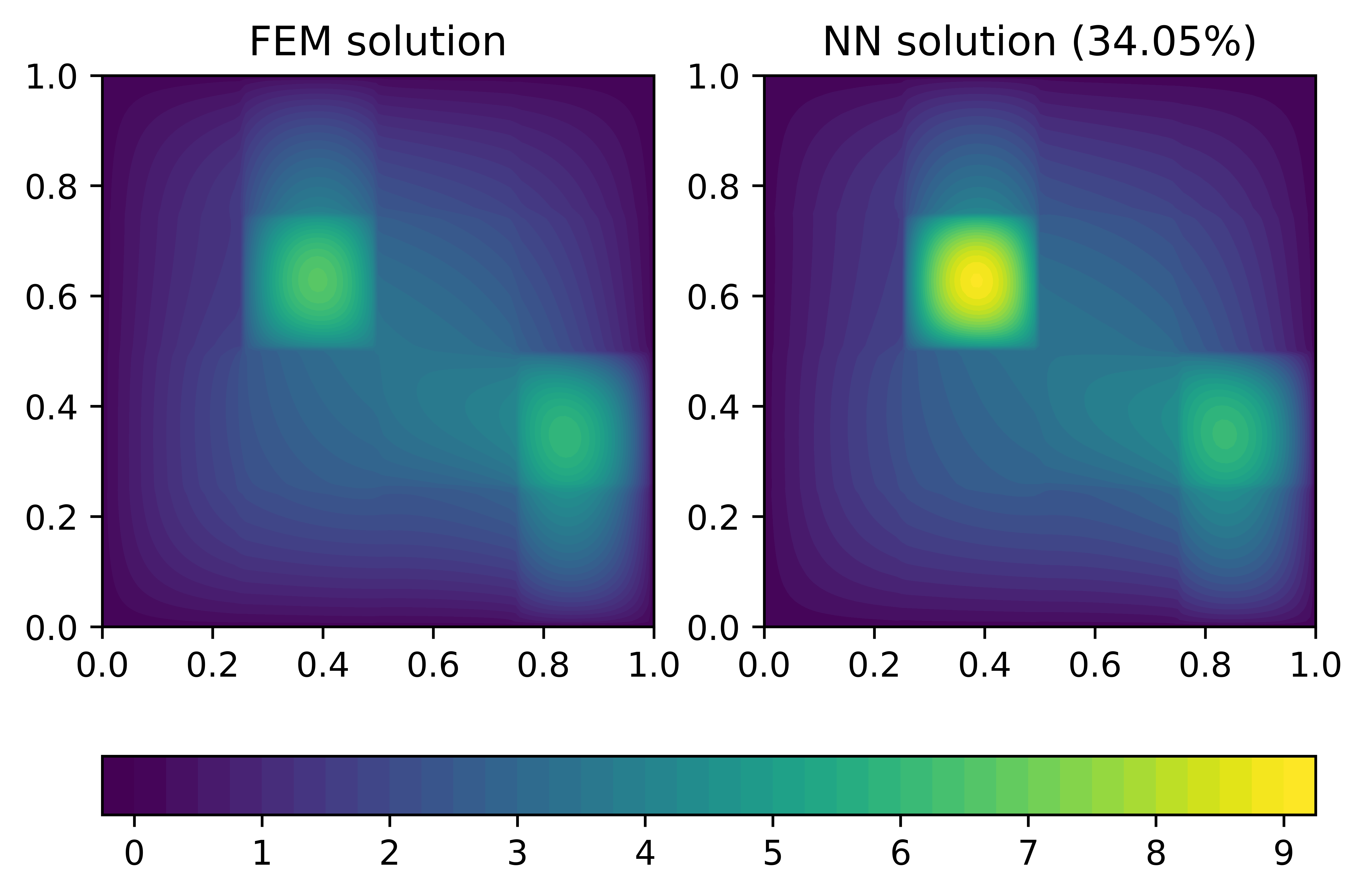}
  \end{tabular}
  \caption{Comparison of the ground truth solution and the one predicted by the NN for an average (left) and a poor performing (right) for $p=16$ and $\mu=10^{-1}$ (top) and $\mu=10^{-3}$ (bottom) for test-case \textbf{[T2]}. The percentage in brackets represents the relative test error for this particular sample.}\label{fig:SquaresExperiment}
\end{figure}

In Figure \ref{fig:SquaresExperiment}, we show samples from the test set for different values of $\mu$. Here we always depict one average performing test-case and one with poor performance. These figures offer a potential explanation of why the scaling with $p$ is qualitatively different for different values of $\mu$. This seems to be because for lower $\mu$ the effect of the individual parameters on the solution seems to be much more local than for higher $\mu$. This appears to lead to a higher intrinsic dimensionality of the problem.

\subsubsection*{\textbf{[T3]} Cookies with Fixed and Variable Radii}

We start with one experiment where the radii of the cookies are fixed to $0.8/(2s)$:
\begin{table}[H]
\renewcommand{\arraystretch}{1.5}
\vspace{-0.5cm}
\small
\begin{tabularx}{\textwidth}{|l|l|l|l|l|l|}
  \hline
  $s$ & 2 & 3 & 4 & 5 & 6 \\ \hline
    \textbf{Parameter dimension $p $} & 4 & 9  & 16 & 25 & 36\\ \hline
  \textbf{Mean relative test error}  & 0.40 \% & 0.41 \% & 0.59 \% & 0.83 \% & 1.10 \%  \\ \hline
  \end{tabularx}
   \caption{Mean relative test error for test-case \textbf{[T3-F]} and different parameter dimensions $p=s^2$ with shift $\mu=10^{-4}$ and radius $0.8/(2s)$.}
     \label{tab:CookiesFixed}
\end{table}
Moreover, we find for the sets of cookies with variable radii $\mathcal{A}^{\mathrm{cvr}}(p,\mu)$ the following mean relative test errors:
\begin{table}[H]
\renewcommand{\arraystretch}{1.5}
\vspace{-0.5cm}
\small
\begin{tabularx}{\linewidth}{|l|l|l|l|l|}
  \hline
  $s$ & 2 & 3 & 4 & 5  \\ \hline
    \textbf{Parameter dimension $p $} & 8 & 18  & 32 & 50 \\ \hline
  \textbf{Mean relative test error }$(\mu=10^{-1})$  & 3.30 \% & 5.44 \% & 7.81 \% & 9.09 \%   \\ \hline
    \textbf{Mean relative test error} ($\mu=10^{-4}$)  & 6.07 \% & 9.81 \% & 12.64 \% &  14.23 \%  \\ \hline
  \end{tabularx}
   \caption{Mean relative test error for test-case \textbf{[T3-V]}  and different parameter dimensions $p=2s^2$.}
     \label{tab:CookiesVarI}
\end{table}

\subsubsection*{\textbf{[T4]} Clipped Polynomials} 

For the set $\mathcal{A}^{\mathrm{cp}}(p,10^{-1})$, we obtain the following mean relative test errors when varying $p$.
\begin{table}[H]
\renewcommand{\arraystretch}{1.5}
\vspace{-0.5cm}
\small
\begin{tabularx}{\textwidth}{|l|l|l|l|l|l|}
  \hline
    \textbf{Polynomial Degree $k$} &  2 & 3 & 5 & 8 & 12 \\ \hline
    \textbf{Parameter dimension $p$} & 6 & 10 & 21 & 45 & 91 \\ \hline
    \textbf{Mean relative test error} & 1.71 \% & 2.58 \% & 3.86 \% & 6.32 \% & 7.58 \% \\ \hline
  \end{tabularx}
     \caption{Mean relative test error for test-case \textbf{[T4]} with clipping value $\mu=10^{-1}$ and different parameter dimensions $p$.}
    \label{tab:ClippedPol}
\end{table}

\begin{figure}[p]
\centering
\begin{minipage}[t]{0.45\textwidth}
    \centering
    \includegraphics[width = \textwidth]{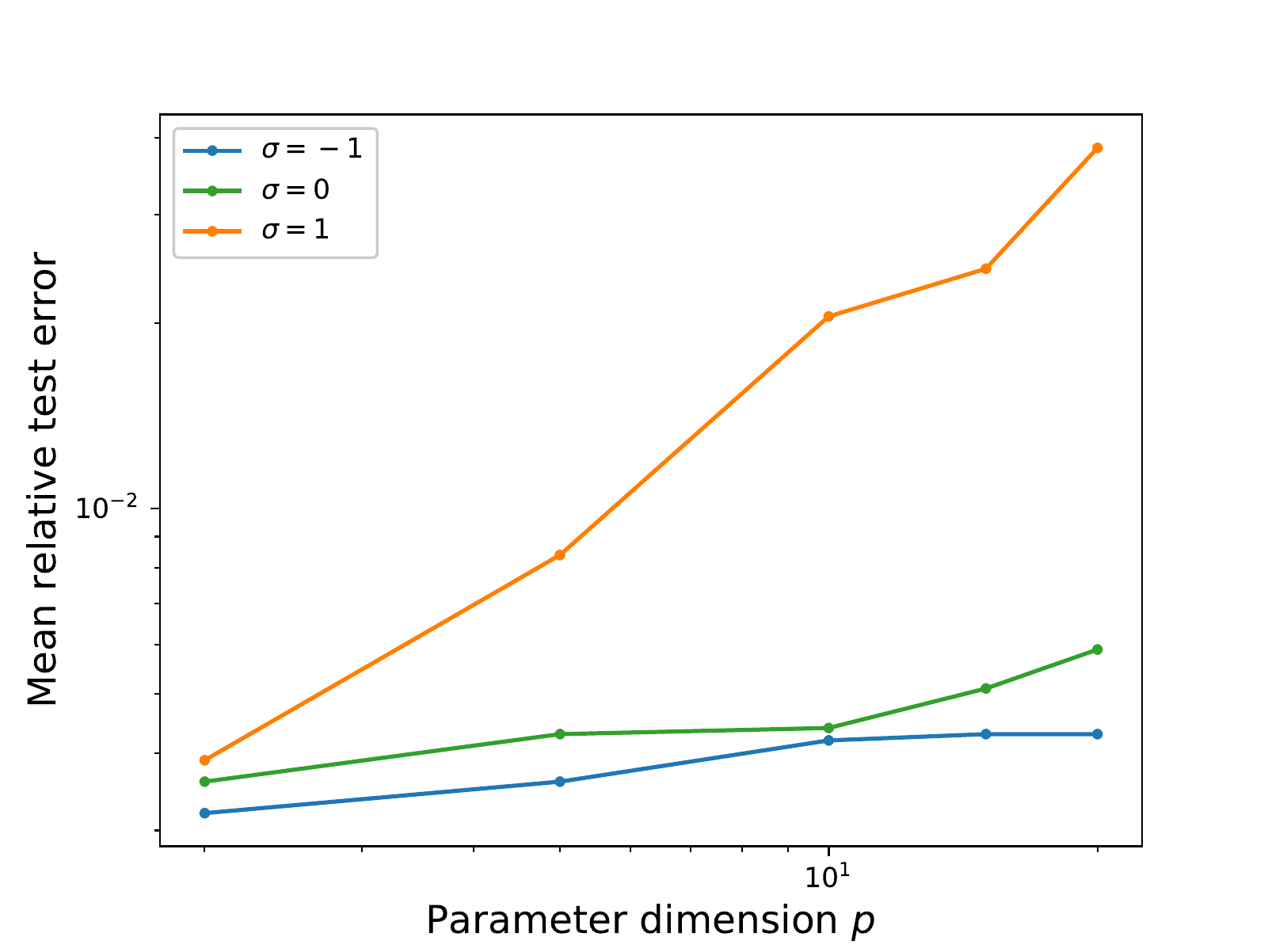}
     \caption{Plot of the mean relative test error for the sets of test-case \textbf{[T1]} for different values $\sigma$. The horizontal axis follows the dimension of the parameter space $p$ the mean relative test error is shown on the vertical axis. Both axes use a logarithmic scale.} 
      \label{fig:trigonometric}
    \end{minipage} \quad
\begin{minipage}[t]{0.45\textwidth}
\centering
  \includegraphics[width = \textwidth]{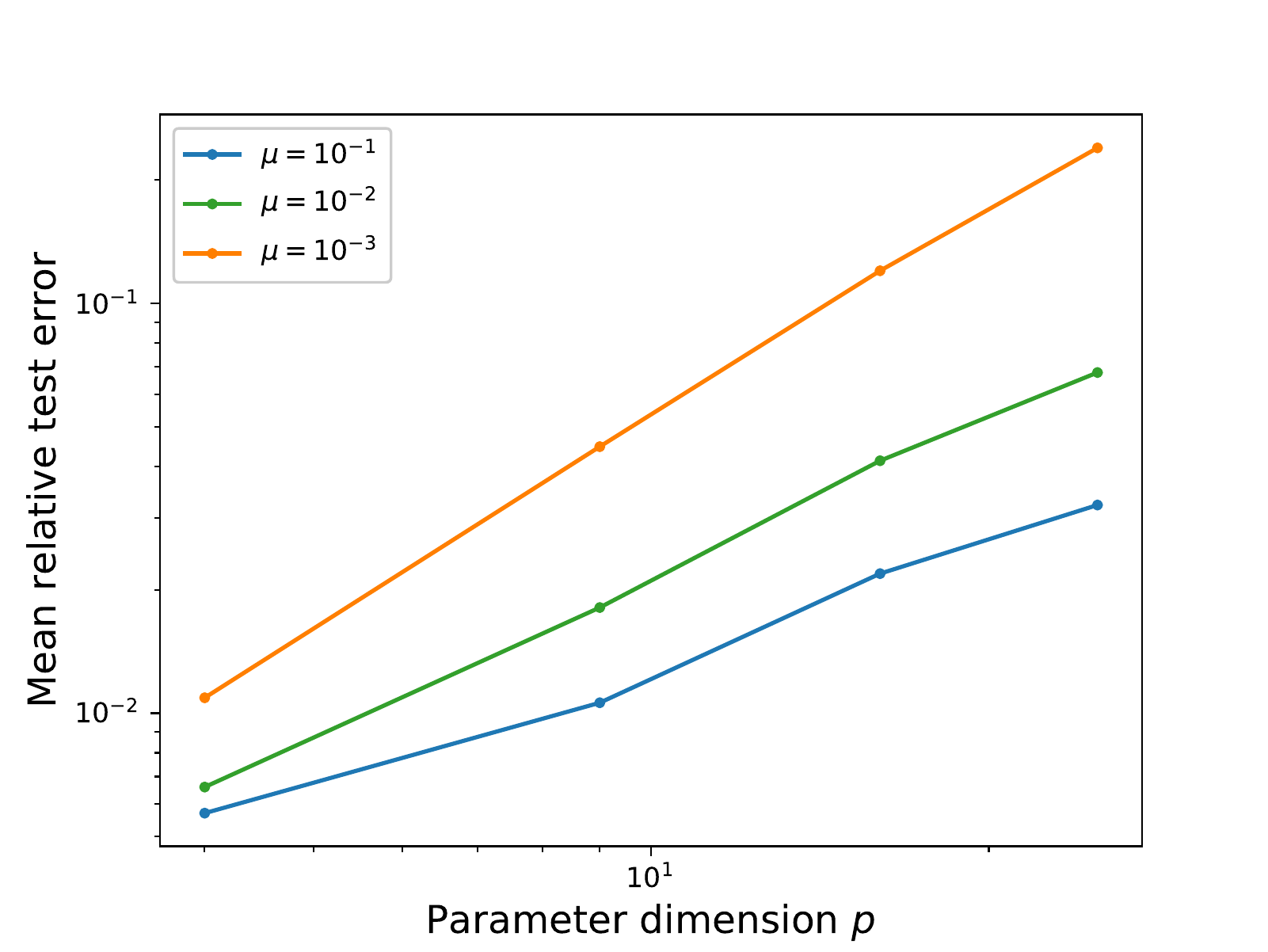}
  \caption{Plot of the mean relative test error for the sets of test-case \textbf{[T2]} for different values of $\mu$ with $p$ on the horizontal axis and the mean relative test error on the vertical axis. Both axes use a logarithmic scale.} 
  \label{fig:chessboard}
\end{minipage} 

\medskip
\begin{minipage}[t]{0.45\textwidth}
    \centering
    \includegraphics[width = \textwidth]{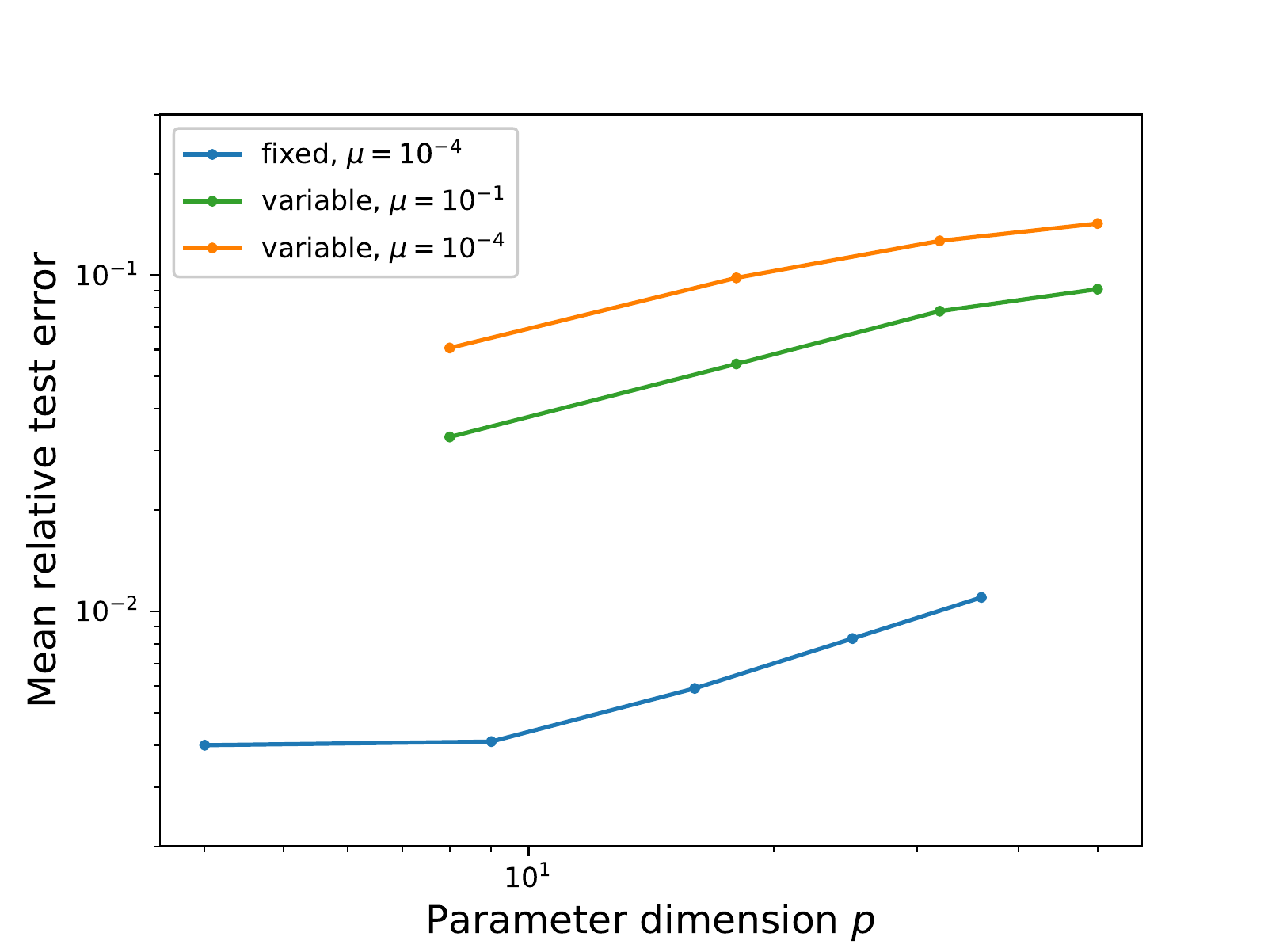}
    \caption{Plot of the mean relative test error for the sets of test-case \textbf{[T3]} with $p$ on the horizontal axis and the error on the vertical axis. Both axes use a logarithmic scale.}
    \label{fig:Cookies}
\end{minipage}\quad
\begin{minipage}[t]{0.45\textwidth}
    \centering
     \includegraphics[width = \textwidth]{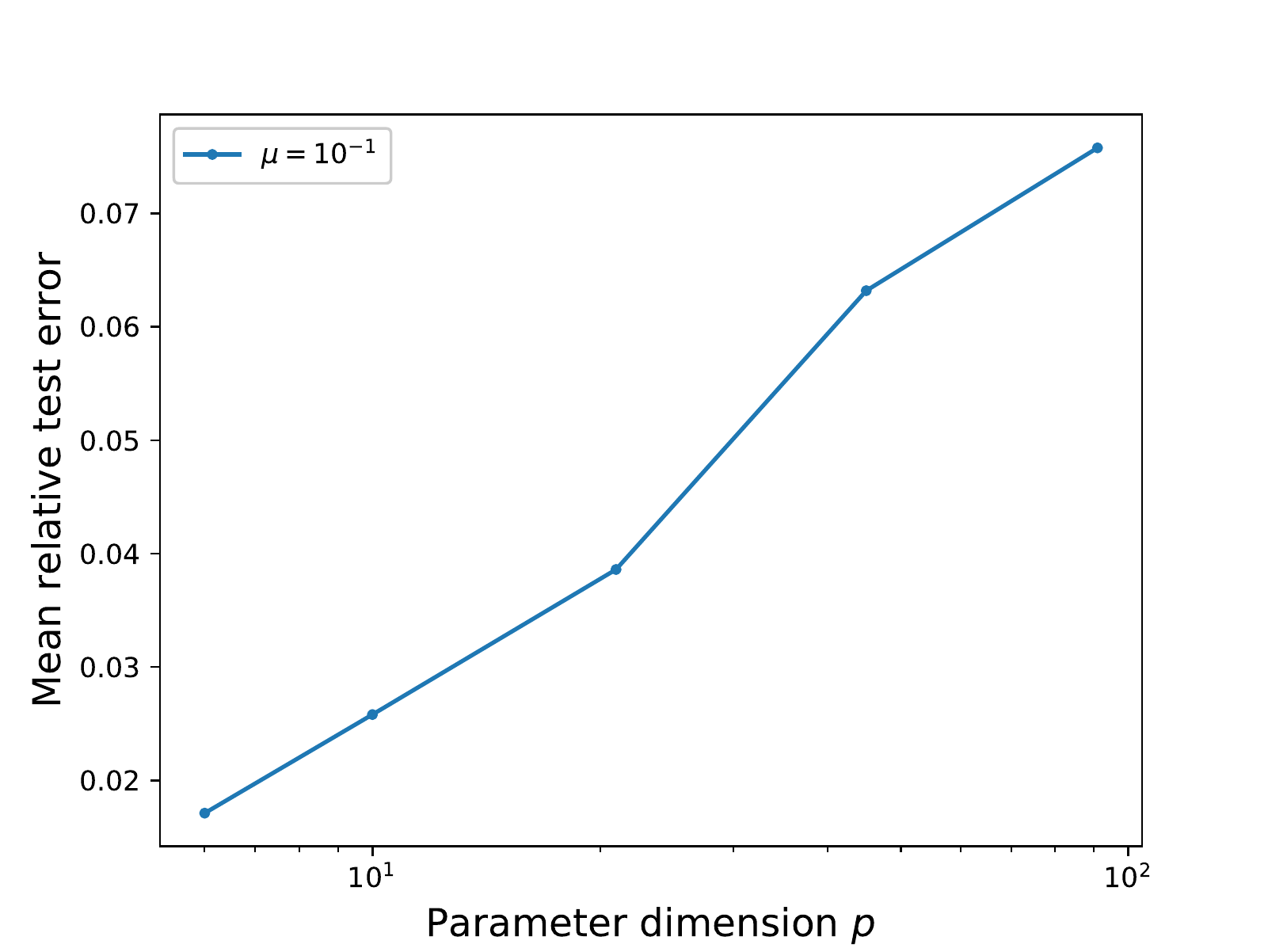}
    \caption{Plot of the mean relative test error for the sets of test-case \textbf{[T4]} with $p$ on the horizontal axis and the error on the vertical axis. Only the horizontal axis is scaled logarithmically.  }
    \label{fig:clippedPoly}
    
    \ \\ \ 
    
    \end{minipage}
\end{figure}

\subsection{Evaluation and Interpretation of Experiments}\label{subsec:Eval}

We make the following observations about the numerical results of Section \ref{subsec:Results}.

\begin{enumerate}
    \item[\textbf{[O1]}] Our test-cases show that the error rate achieved by NN approximations for varying parameter sizes \emph{differs strongly and qualitatively between different test-cases}. In Figures \ref{fig:trigonometric}, \ref{fig:chessboard}, \ref{fig:Cookies}, and \ref{fig:clippedPoly} we depict the different scaling behaviors of the test-cases \textbf{[T1]}, \textbf{[T2]}, \textbf{[T3]}, and \textbf{[T4]}. For \textbf{[T1]} and $\sigma = -1$, \emph{the error appears to be almost independent from $p$} for $p \to \infty$. In contrast to that, we observe for $\sigma = 1$ a linear scaling in the loglog plot implying a polynomial dependence of the error on $p$.
    
    For test-case \textbf{[T2]}, we observe that the error scales linearly in the loglog scale of Figure \ref{fig:chessboard}. We conclude that for $\mathcal{A}^{\mathrm{cb}}(p,\mu),$ \emph{the error scales polynomially with $p$}. 
    
    The errors of the test-cases associated with \textbf{[T3]} seem to scale linearly with $p$ in the loglog scale depicted in Figure  \ref{fig:Cookies}. This implies that for \textbf{[T3]} \emph{the error scales polynomially in $p$ with the same exponent.}  
    
    The semilog plot of Figure \ref{fig:clippedPoly} shows that for test-case \textbf{[T4]} with the  sets $\mathcal{A}^{\mathrm{cp}}(p,10^{-4}),$ \emph{the growth of the error is logarithmic in $p$}. 

    In total, we observed scaling behaviors of $\mathcal{O}(1), \mathcal{O}(\log(p))$ and $\mathcal{O}(p^k)$ for $k>0$ and for $p \to \infty$. Notably, none of the test-cases exhibited an exponential dependence of the error on $p$.

    \item[\textbf{[O2]}] The choice of the hyper-parameters $\sigma$ and $\mu$ in the test-cases \textbf{[T1]}, \textbf{[T2]}, \textbf{[T3]} \emph{influences the scaling behavior according to its effect on the complexity of the parameterized diffusion coefficient set}. 

    Weighting the parameters using the scaling parameter $\sigma$ should, in principle, \emph{simplify the parametric problem for smaller values of $\sigma$}. This is precisely, what we observe in Table \ref{tab:TrigonPol} and Figure \ref{fig:trigonometric}.

    The influence of the shift $\mu$ is of a somewhat different type. \emph{Higher values of $\mu$ make the underlying problem more elliptic}. This can be seen in Figure \ref{fig:SquaresExperiment}: For a small value of $\mu$, the impacts of the individual values on the chessboard-pieces on the solution appear to be almost completely decoupled. On the other hand, in the more elliptic case, the solution appears more smoothed out, and therefore each parameter value also influences the solution more globally. This implies a stronger coupling of the parameters and at least intuitively indicates a \emph{reduced intrinsic dimensionality for higher values of $\mu$}. 

    Accordingly, we see in Table \ref{tab:SquareI} and Figure \ref{fig:chessboard} that the parameter $\mu$ influences the scaling behavior of the method with $p$. Indeed, the error scales as $\mathcal{O}(p^k)$, where the exponent in the polynomial dependence on $p$ depends on $\mu$.
    
    Concluding, we can see that the approximation of the DPtSM by NNs appears to be very sensitive to these parameters, as we observe in Table \ref{tab:SquareI} and Figure \ref{fig:chessboard} as well as in Table \ref{tab:CookiesVarI} and Figure \ref{fig:Cookies}. 

    \item[\textbf{[O3]}] We observe \emph{no fundamentally worse scaling behavior for non-affinely parametrized test-cases} compared to test-cases with an affine parameterization. In test-case \textbf{[T3]}, we do observe that the non-linearly parametrized problem appears to be more challenging overall, while the scaling behavior is the same as for the affinely parametrized problem. 
    In test-case \textbf{[T4]}, which is the test-case with the highest number of parameters $p$, we observe only a very mild (in fact logarithmic) dependence of the error on $p$.
\end{enumerate}

From these observations we draw the following conclusions for our hypotheses:

\paragraph{Hypothesis \textbf{[H1]}}

In observation \textbf{[O1]}, we saw that over a wide variety of test-cases multiple types of scaling of the error with the dimension of the parameter space could be observed. None of them admit an exponential scaling. In fact, the behavior of the errors seems to be determined by an intrinsic complexity of the problems. 
\paragraph{{Hypothesis} \textbf{[H2]}}

Comparing performance both between test-cases (observation \textbf{[O1]}) and within test-cases (observation \textbf{[O2]}), leads us to conclude that there exist strong differences in the performance of learning the DPtSM. 
For various test-cases, using NNs with precisely the same architecture, we observed (see \textbf{[O2]}) considerably different scaling behaviors of the test-cases \textbf{[T1]}-\textbf{[T4]} which have the error scale polynomially, logarithmically and being constant with changing parameter dimension $p$ (described in \textbf{[O1]}).  
According to \textbf{[O2]}, the overall level of the errors and the type of scaling for increasing $p$ follows the semi-ordering of complexities of test-cases in the sense that more complex parametrized sets yield higher errors whereas simpler sets or spaces with intuitively lower intrinsic dimensionality yield smaller errors (test-cases \textbf{[T1]} and \textbf{[T2]}). 

Therefore, we conclude that the approximation theoretical intrinsic dimension of the parametric problem is a main factor in determining the hardness of learning the DPtSM.

\paragraph{{Hypothesis} \textbf{[H3]}}

In support of \textbf{[H3]}, we found no fundamental difference of the performance of the NN model for non-affinely parametrized problems (see \textbf{[O3]}).

\medskip

In conclusion, we found support for all the hypotheses \textbf{[H1]}-\textbf{[H3]}. We consider this result a validation of the importance of approximation-theoretical results for practical learning problems, especially in the application of deep learning to problems of numerical analysis.

It is clear that the results presented in this work only analyze the sensitivity of the performance of the learned DPtSM corresponding to the semi-ordering of complexities. For future work, it would be interesting to identify alternative and more quantitative notions of complexities and test the sensitivity of the learned method with regards to those.

\section*{Acknowledgements}

M. Geist and M. Raslan would like to thank Philipp Trunschke for fruitful discussions on the topic. This work was made possible by the computational resources provided by the Institute of Mathematics of the TU Berlin. 
G. Kutyniok acknowledges partial support by the Bundesministerium
 f\"ur Bildung und Forschung (BMBF) through the Berlin Institute for the
Foundations of Learning and Data (BIFOLD), Project AP4, RTG DAEDALUS (RTG 2433),
Projects P1, P3, and P8, RTG BIOQIC (RTG 2260), Projects P4 and P9, and by the
Berlin Mathematics Research Center MATH+, Projects EF1-1 and EF1-4.

\small
\bibliographystyle{abbrv}
\bibliography{references}
\normalsize

\appendix
\section{Elimination of Obfuscating Phenomena}\label{sec:EliminationCauses}

Below we describe the measures taken to enable comparability between test-cases. 

\subsection{Fixing the Architecture}\label{subsec:FixArch}

In all our experiments the network architecture was kept almost completely fixed, only varying the dimension of the input layer. Our choice of architecture was made on the basis of preliminary experiments with the goal of developing a network structure that performs well on all datasets and in particular displays good optimization behavior independent of the test-case as showcased in Appendix \ref{subsec:OptBehav}. This was done to ensure comparability across all test-cases and parameter choices, allowing us to isolate the influence of the parametrization and the dimension of the parameter space. We emphasize that more sophisticated architectures and the usage of tools like weight regularization or learning rate decay in general enable better performance on individual datasets. However, in our case, they would only obfuscate the approximation-theoretical effect that we are seeking to identify.        

\subsection{Influence of the Size of the Training Set} \label{subsec:Sample}

Throughout this paper all training was conducted with a fixed number of 20000 samples. Since it is clear that a larger training set will generally yield better results, this trend may affect different test-cases to various degrees. To guarantee that the effect of the choice of the number of training samples is uniform across cases, we chose the number of samples in the following way:
We trained the same NN architecture as described in Subsection \ref{sec:SetupNN} for different parameter constellations with training sets ranging from 10000 to 20000 samples. The results are depicted in Table \ref{tab:SampleVar}. The table also includes the coefficient of determination $R^{2}$ (see \cite[p. 601]{R2}) for each individual dataset resulting from fitting a simple linear regression to the set of sample size and test error pairs.  

\begin{centering}
\renewcommand{\arraystretch}{1.5}
\vspace{-0.2cm}
\small
\begin{tabularx}{\linewidth}{|l|l|l|l|l|l|l|}
  \hline
  	\diagbox[width=5.1cm,height=1.8cm]{\textbf{Test-case}}{\textbf{Size of training set}}  & 20000 & 17500 & 15000 & 12500 & 10000 & $R^{2}$ \\ \hline
 {\textbf{[T1]}} ($\mu=1,$ $\sigma=0,$ $p=20$)  & 0.59 \% & 0.61 \% & 0.64 \% &  0.70 \%   &  0.76 \% &  0.95 \\ \hline
  {\textbf{[T2]}} ($\mu=10^{-1},$ $p=9$)  & 1.06 \% & 1.29 \% & 1.49 \% &  1.81 \%   &  2.18 \% & 0.98 \\ \hline
{\textbf{[T2]}} ($\mu=10^{-2},$ $p=9$)  & 1.81 \% & 1.94 \% & 2.58 \% &  2.98 \%   &  4.26 \% & 0.91 \\ \hline 
{\textbf{[T2]}} ($\mu=10^{-3},$ $p=9$)  & 4.47 \% & 5.31 \% & 6.23 \% &  7.78 \%   &  9.24 \% & 0.98 \\ \hline
 {\textbf{[T3-F]}} ($\mu=10^{-4},$ $p=25$)  & 0.83 \% & 0.85 \% & 0.88 \% &  0.91 \%   &  0.96 \% & 0.97 \\ \hline
  {\textbf{[T3-V]}} ($\mu=10^{-1},$ $p=18$)  & 5.44 \% & 5.60 \% & 5.83 \% & 6.16\%   &  6.56 \% & 0.97 \\ \hline 
  {\textbf{[T3-V]}} ($\mu=10^{-4},$ $p=18$)  & 9.81 \% & 9.98 \% & 10.18 \% & 10.61\%   &  11.06 \% & 0.95 \\ \hline 
  {\textbf{[T4]}} ($\mu=10^{-1},$ $p=21$)  & 3.86 \% & 4.17 \% & 5.06 \% &  5.50 \%   &  6.46 \% & 0.98 \\ \hline
   \caption{Mean relative test error as well as the corresponding $R^{2}$ coefficient from a simple linear regression for varying sizes of the training set and all previously considered setups.}
     \label{tab:SampleVar}
  \end{tabularx}
\end{centering}

This analysis shows that with $R^{2}$ values ranging from 0.91 to 0.98 the relation between the number of samples and the achieved accuracy is almost perfectly linear. Assuming this relation extrapolates to the other parameter dimension $p$, this implies that our results in Section \ref{sec:NumericalRes} can be considered independent of the number of samples chosen. It should, however, be noted, that this linear dependence can only be observed in a reasonable range of training set sizes. In particular, the experiments revealed a lower bound on the number of samples needed to stably train our NN architecture. While in our case this bound can be observed in the range of 1000 to 5000 samples depending on the considered test-case, other NN setups may be able to effectively train with even lower sample counts.  

\subsection{A Posteriori Analysis of Convergence Behavior}\label{subsec:OptBehav}

\begin{figure}[htb]
\centering
\begin{minipage}[t]{0.45\textwidth}
    \centering
    \includegraphics[width = \textwidth]{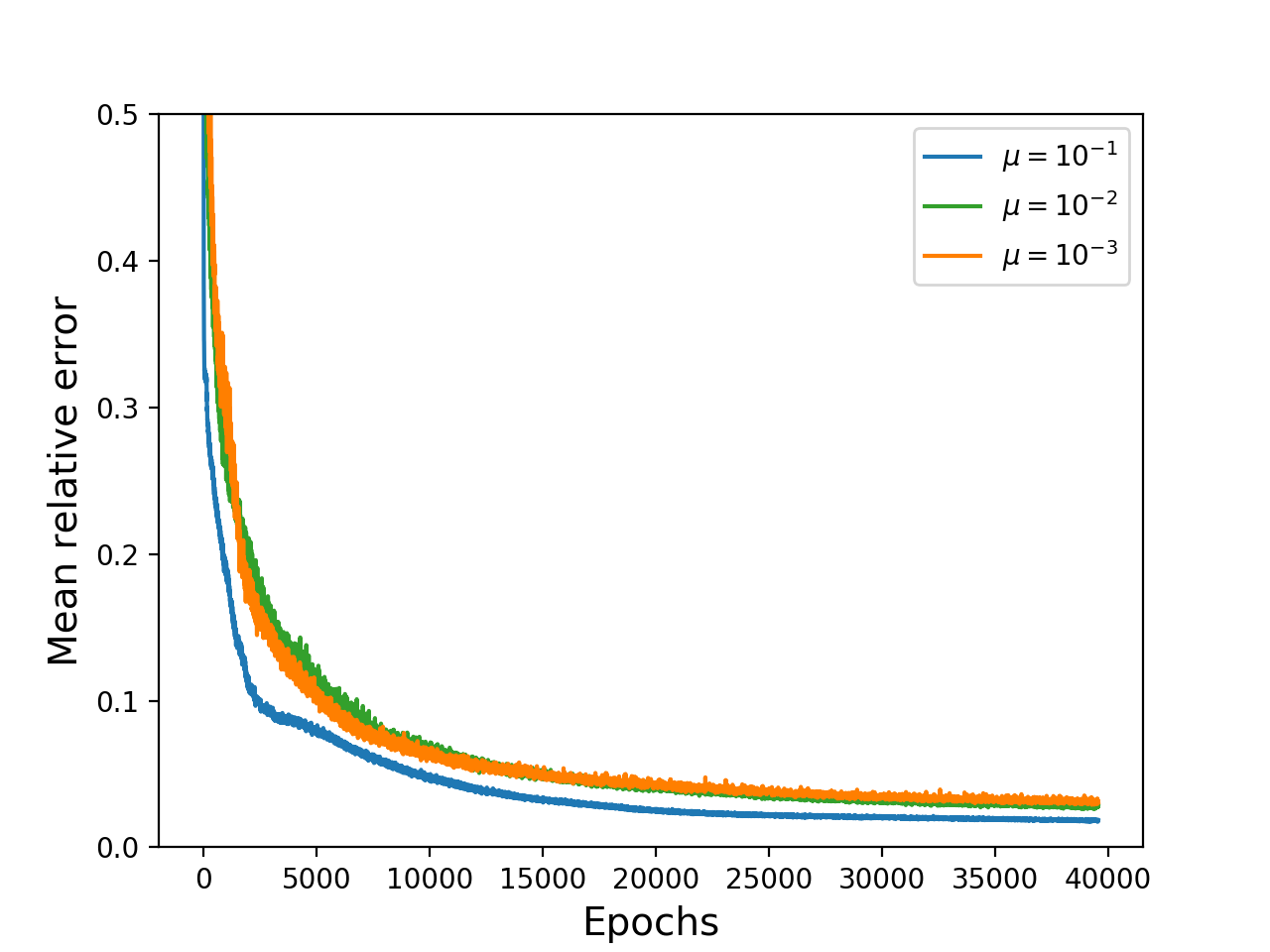}
     \caption{Plot of the mean relative training error for \textbf{[T2]} with $p=25$ and different shifts $\mu$.} 
      \label{fig:conv_squares}
    \end{minipage} \quad
\begin{minipage}[t]{0.45\textwidth}
\centering
  \includegraphics[width = \textwidth]{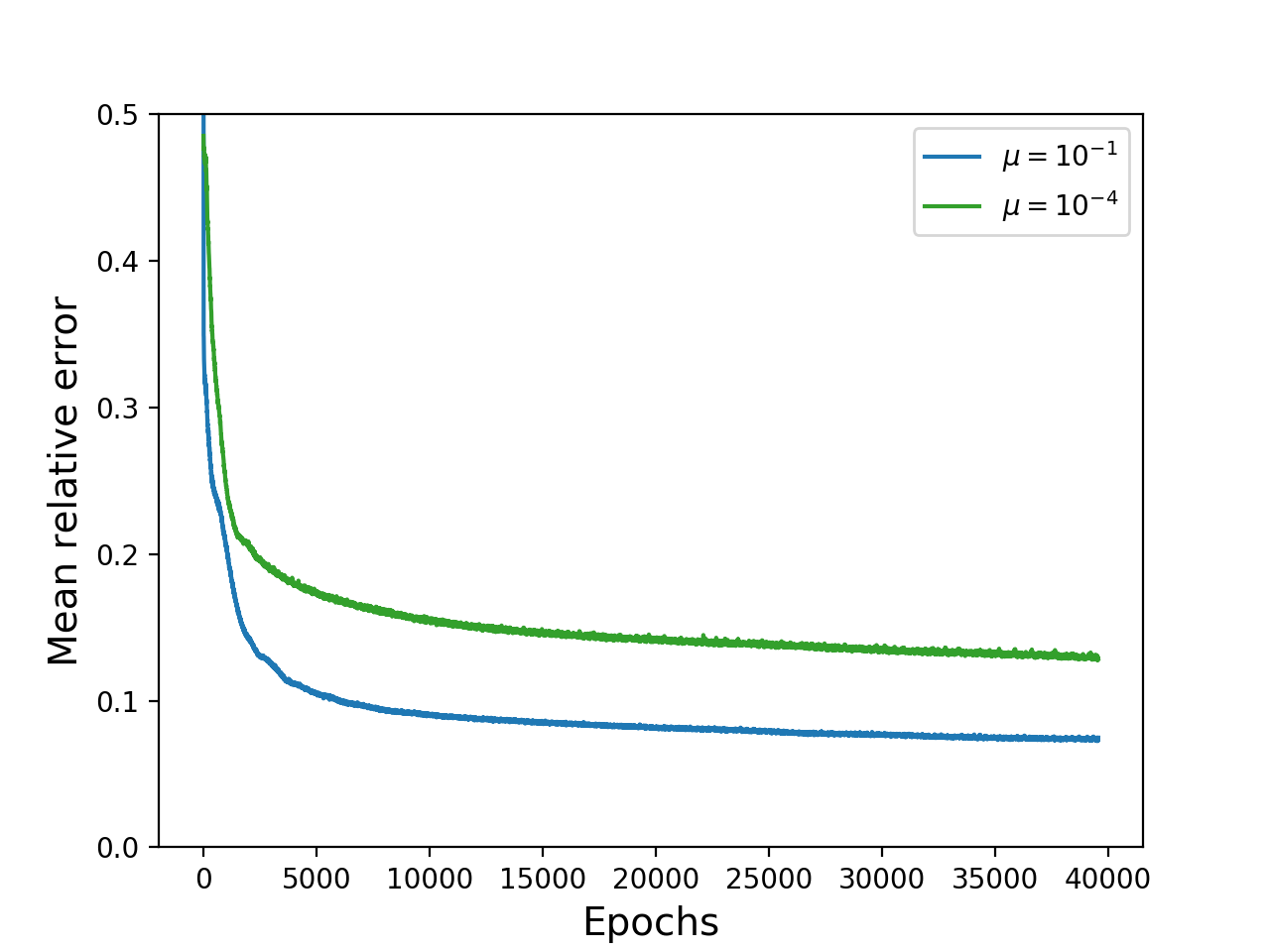}
  \caption{Plot of the mean relative training error for \textbf{[T3-V]} with $p=50$ and different shifts $\mu$.} 
  \label{fig:conv:cookies}
\end{minipage} 

\end{figure}

Similarly to the architecture, the hyper-parameters of the optimization method were also kept fixed across all datasets and training runs. This measure, however, only eliminates the effect of the architecture on the optimization method and does not address any obfuscating effect that the choice of test-cases may have. To analyze if such an effect is present, we check the convergence on our two hardest test-cases \textbf{[T2]} and \textbf{[T3-V]} for the largest parameter dimension $p$ considered. The results are depicted in Figure \ref{fig:conv_squares} and \ref{fig:conv:cookies}, respectively. We see that even for small shifts $\mu$, i.e., the most difficult problem settings, the error on the training set converges smoothly. This behavior can also be witnessed on all other test-cases.

\begin{figure}[htb]
\centering
\begin{minipage}[t]{0.45\textwidth}
    \centering
    \includegraphics[width = \textwidth]{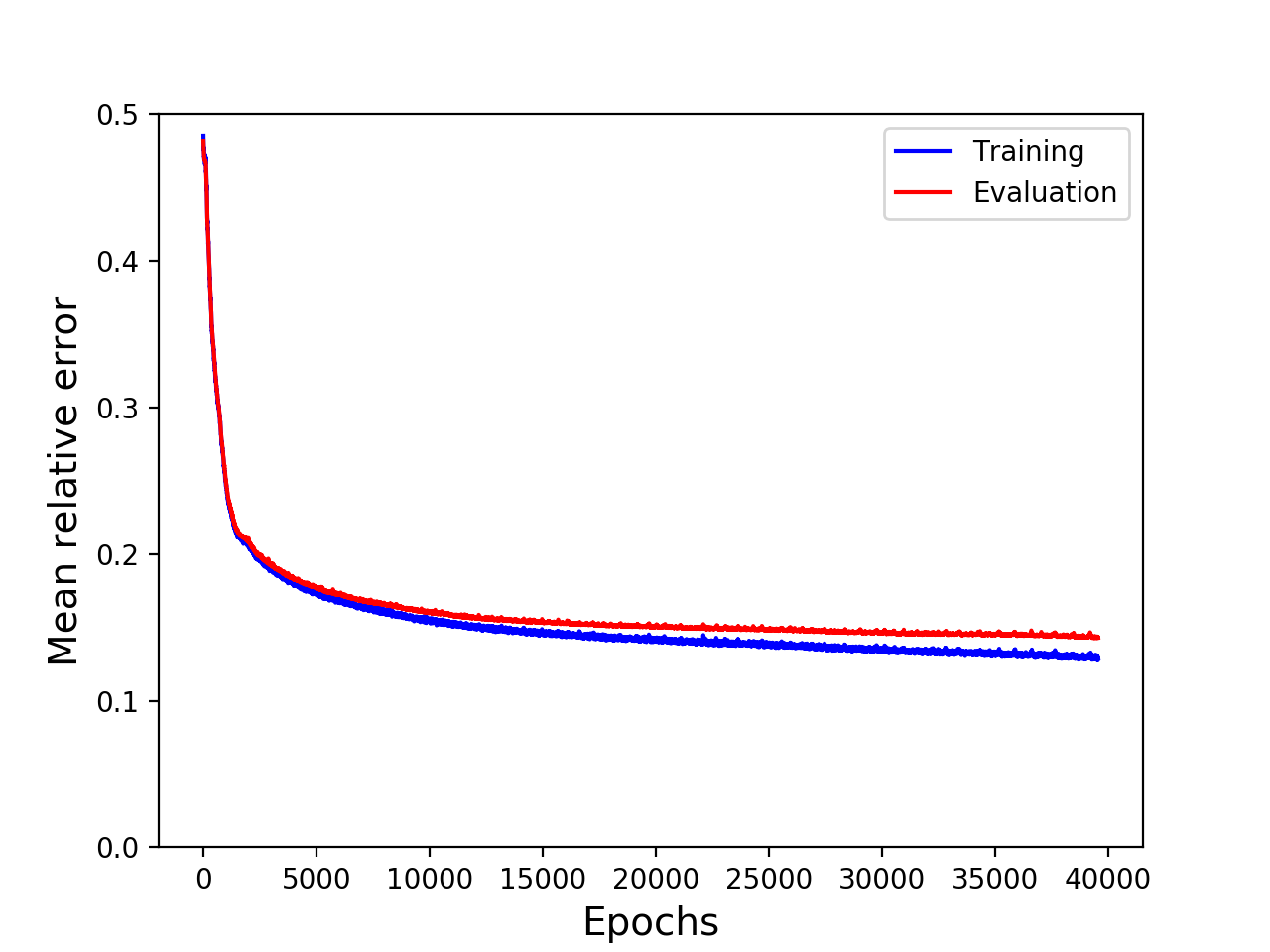}
     \caption{Plot of the mean relative training and test error for \textbf{[T3-V]} with $p=50$ and $\mu=10^{-4}$.} 
      \label{fig:conv_cookies_eval}
    \end{minipage} \quad
\begin{minipage}[t]{0.45\textwidth}
\centering
  \includegraphics[width = \textwidth]{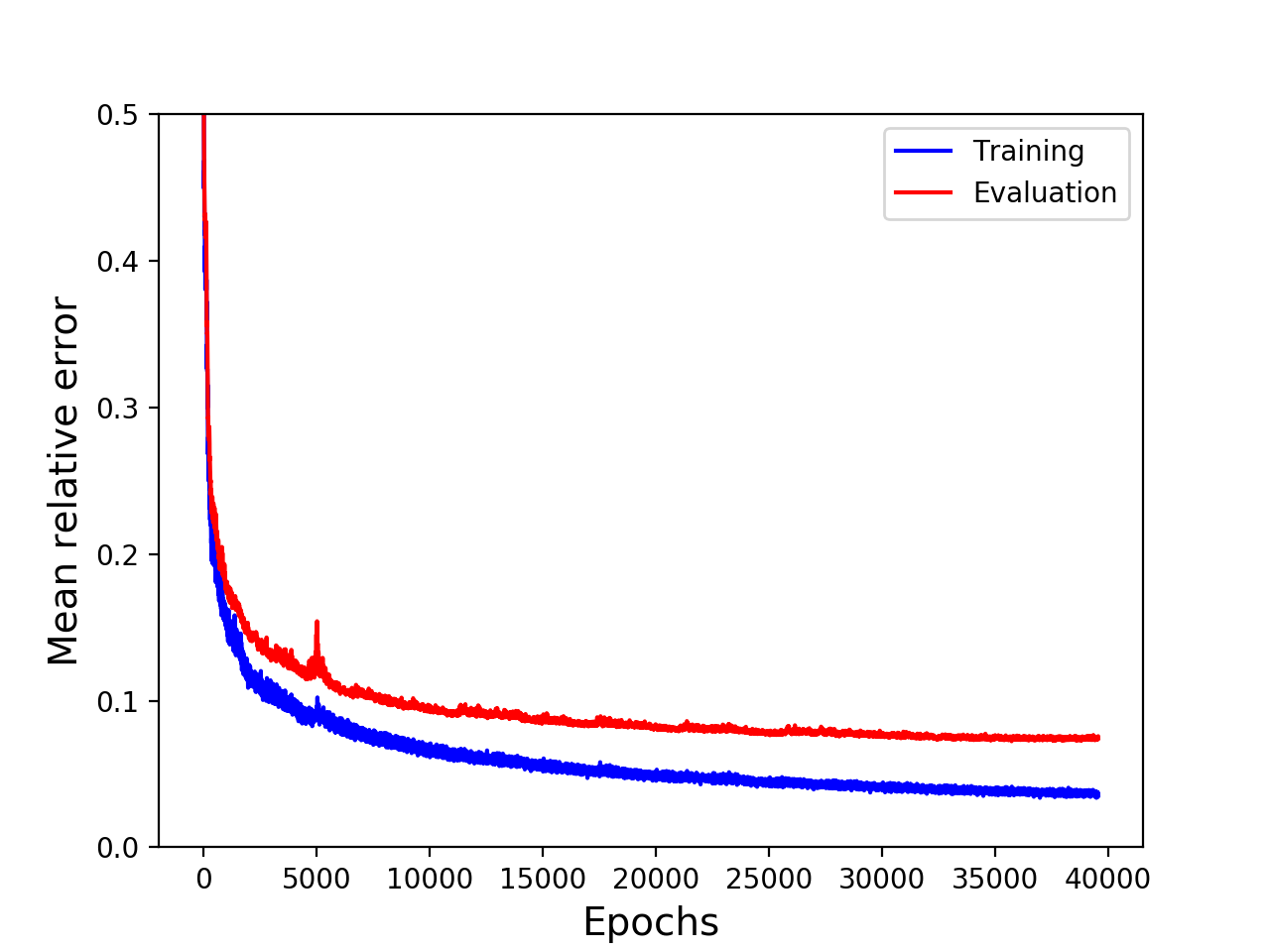}
  \caption{Plot of the mean relative training and test error for \textbf{[T4]} with $p=91$ and $\mu=10^{-1}$.} 
  \label{fig:conv_poly_eval}
\end{minipage} 

\end{figure}

Another possible pitfall of our optimization procedure would be the occurrence of overfitting. In particular, this would render our attained accuracy levels invalid as we trained for a fixed number of epochs. However, this did not occur in any of our tests. We exemplarily showcase the convergence plot of the training and test error for the hardest parameter choices of \textbf{[T3-V]} and \textbf{[T4]} in Figure \ref{fig:conv_cookies_eval} and \ref{fig:conv_poly_eval} respectively. Similar behavior can also be observed on all other datasets.  
 
\end{document}